\documentclass[hidelinks,onefignum,onetabnum]{siamart251216}

\usepackage{graphicx}
\usepackage{algorithmicx}
\usepackage[noend]{algpseudocode}
\usepackage{bm}
\usepackage{amsfonts}
\usepackage{float}
\usepackage{epstopdf}
\usepackage{booktabs}
\usepackage{rotating}

\ifpdf
\DeclareGraphicsExtensions{.pdf,.png,.jpg,.eps}
\else
\DeclareGraphicsExtensions{.eps}
\fi

\newcommand{\R}{\mathbb{R}}
\newcommand{\tr}{^{\sf T}}
\newcommand{\m}[1]{{\bf{#1}}}
\newcommand{\g}[1]{\bm #1}

\newsiamremark{remark}{Remark}
\newsiamthm{assum}{Assumption}
\crefname{assum}{Assumption}{Assumptions}
\definecolor{noteyellow}{rgb}{0.58,0.40,0.00}
\newcommand{\notelabel}{{\color{noteyellow}\textbf{Note.}} }

\newcommand{\redfix}[1]{#1}
\newcommand{\bluefix}[1]{}
\newcommand{\newredfix}[1]{#1}
\newcommand{\yc}[1]{{\color{noteyellow}#1}}
\newcommand{\Cstd}{{\mathcal{C}_{\mathrm{std}}}}
\newcommand{\Cmod}{{\mathcal{C}_{\mathrm{mod}}}}
\newcommand{\Chor}{{\mathcal{C}_{\mathrm{hor}}}}
\newcommand{\Cmodraw}{{\widetilde{\mathcal{C}}_{\mathrm{mod}}}}

\headers{Optimal Last-Iterate Rates for PDSG}
{H. Wu, C. Yin, Y. Dai, and L. Wang}
\title{Optimal Last-Iterate Rates for the Projected Decentralized Stochastic Subgradient Method \thanks{This work was supported by the National Natural Science Foundation of China (Grant Number 11971231, 92473208).}}
\author{
	Hao Wu\thanks{Academy of Mathematics and Systems Science, Chinese Academy of Sciences, Beijing 100190, China
		(\email{wuhoo104@amss.ac.cn}, \email{dyh@lsec.cc.ac.cn}).}
	\and Chao Yin\thanks{School of Mathematics, Hohai University, Nanjing 211100,
		Jiangsu, China (\email{yinchao@hhu.edu.cn}).}
	\and Yuhong Dai\footnotemark[2]\hspace{0.42em}\thanks{School of Mathematical Sciences,
		University of Chinese Academy of Sciences, Beijing 100049, China.}
	\and Liping Wang\thanks{Corresponding author. School of Mathematics,
		Nanjing University of Aeronautics and Astronautics, Nanjing 211106,
		Jiangsu, China (\email{wlpmath@nuaa.edu.cn}).}
}

\ifpdf
\hypersetup{
	pdftitle={Optimal Last-Iterate Rates for the Projected Decentralized Stochastic Subgradient Method},
	pdfauthor={Hao Wu, Chao Yin, Yuhong Dai, Liping Wang}
}
\fi

\begin{document}
	\maketitle
	
	\begin{abstract}
We study the projected decentralized stochastic subgradient (PDSG) method for minimizing the average of convex, possibly nonsmooth local
objectives over a common closed convex constraint set and a fixed
connected undirected network. We focus on the finite-time expected suboptimality of the last
local iterate generated at each node, which is the actual output locally
available without time averaging or network averaging. Under the standard stepsizes $\alpha_t=\alpha/\sqrt{t}$ and $\alpha_t=2/(\mu t)$,
{we establish} $\mathcal{O}(\log T/\sqrt{T})$ and
$\mathcal{O}(\log T/T)$ last-local-iterate rates for PDSG when each local objective is convex and the aggregate objective is strongly convex, respectively.
To eliminate the logarithmic loss, we introduce a general
horizon-dependent stepsize framework, termed cooldown stepsize, that combines a base stepsize
with segmentwise scaling over a multistage partition and includes the dyadic construction (Jain et al, SIAM J. Optim., 2021) as a special case. For adapting to decentralized settings, we propose a long-range variation condition to make the stepsize schedule compatible with network contraction. 
Under this cooldown stepsize, we obtain $O(1/\sqrt{T})$ and $O(1/T)$ last-local-iterate rates for PDSG in both problem classes, matching the optimal centralized dependence on $T$. 
Experiments on a lower-bound instance and constrained hinge-loss problems illustrate the predicted terminal behavior of local iterates and the control of consensus error under the cooldown stepsize.
	\end{abstract}
	
	\begin{keywords}
		decentralized optimization, stochastic subgradient method, last-iterate convergence, convex optimization, cooldown stepsize
	\end{keywords}
	
	\begin{MSCcodes}
		90C15, 90C25, 90C30, 65K05
	\end{MSCcodes}
	\section{Introduction}\label{sec:introduction}
	Decentralized optimization enables a network of nodes to minimize an aggregate objective without a central coordinator. Each node accesses only local information and communicates only with its neighbors. This paradigm arises in resource allocation \cite{fusco2021decentralized}, wireless networks \cite{jeong2022asynchronous}, decentralized machine learning \cite{zhang2022distributed}, power systems \cite{weitenberg2018robust}, and federated learning \cite{lei2025decentralized}.
	We consider $n$ nodes communicating over a fixed, connected, undirected graph $\mathcal{G}=(\mathcal{V},\mathcal{E})$, where $\mathcal{V}=\{1,\ldots,n\}$. {Each node $i$ has access only to its local convex, possibly nonsmooth objective $f_i:\mathbb{R}^p\to\mathbb{R}$, and all nodes {share a common constraint set}.} {They seek a common minimizer of}
	\begin{equation}\label{obj_fun1}
		\min_{\m{z}\in\Omega}\; F(\m{z}):=\frac{1}{n}\sum_{i=1}^n f_i(\m{z}),
	\end{equation}
	{where $\Omega\subseteq\mathbb{R}^p$ is a nonempty, closed, and convex
	set. We assume that the solution set of
	\eqref{obj_fun1} is nonempty, and fix an arbitrary optimizer $\m z^*\in\arg\min_{\m z\in
\Omega}F(\m z)$ and write $F^*=F(\m z^*)$.}
	
	A classical method for \eqref{obj_fun1} is the projected decentralized stochastic subgradient (PDSG) method \cite{nedic2009distributed, Sundhar2010}. Each node maintains a local copy $\m{x}_i^t$ and updates it by
	\begin{align}\label{update_rule}
		\m{x}_i^{t+1}=\mathcal{P}_\Omega\left[\sum_{j=1}^n\tilde W_{ij}\m{x}_j^t-\alpha^t\m{g}_i^t\right],
	\end{align}
	where $\alpha^t>0$ is the stepsize, $\tilde{\m W}=[\tilde W_{ij}]$ is a mixing matrix, $\mathcal{P}_\Omega(\cdot)$ denotes the Euclidean projection onto $\Omega$, and $\m g_i^t$ is a conditionally unbiased estimator of a subgradient of $f_i$ at $\m x_i^t$. In terms of the stacked vectors $\m{x}^t=[\m{x}_1^t;\ldots;\m{x}_n^t]$ and $\m{g}^t=[\m{g}_1^t;\ldots;\m{g}_n^t]$, the update becomes
	$
		\m{x}^{t+1}=\mathcal{P}_{\Omega^n}\left[\m W\m{x}^t-\alpha^t\m{g}^t\right],
	$
	where $\m W=\tilde{\m W}\otimes\m I_p$ and $\Omega^n=\underbrace{\Omega \times \ldots \times \Omega}_{n}$.

	
	{Finite-time guarantees for classical decentralized stochastic subgradient methods are usually formulated in terms of network averages $\bar{\m{x}}^t$, time-averaged iterates, or best iterates. These quantities, however, do not directly characterize the native state generated by the algorithm at an individual node.}
	Indeed, computing $\bar{\m{x}}^t$ at a node requires additional network-wide aggregation, whereas time averaging or best-iterate selection requires  additional variables storage.
    These considerations motivate a finite-time analysis of the last local iterate $\m{x}_i^T$, where $T$ denotes the iteration horizon.

	{Last-iterate convergence is known to be more delicate
	than ergodic convergence in nonsmooth stochastic optimization.}
	In the centralized setting, {standard schedules such as $\alpha/\sqrt{t}$ and $2/(\mu t)$ can incur, }
    in the worst case, an additional logarithmic factor relative to the corresponding lower bounds \cite{hazan2011beyond,shamir2013stochastic,rakhlin2012making,harvey2019tight}. Recent work has developed modified stepsize schedules that recover the optimal rates \cite{jain2021making,zamani2025exact}.
	{Nonergodic convergence has also been investigated in decentralized
	nonsmooth optimization. In particular, Dong et al.~
	\cite{dong2021nonergodic} established
	$\mathcal{O}(\log T/\sqrt{T})$ and
	$\mathcal{O}(\log T/T)$ rates for the last iterate of a decentralized stochastic proximal-gradient method in the convex and strongly convex settings, respectively, whereas their optimal-order
	$\mathcal{O}(1/\sqrt{T})$ and $\mathcal{O}(1/T)$ guarantees are obtained through time averaging. Thus, the existing theory leaves a gap between the rates achievable by time-averaged iterate and by the last iterate itself.}
	
	{Closing this gap for the last local iterate is nontrivial in
	decentralized optimization. In addition to the last-iterate difficulty already present in centralized nonsmooth stochastic optimization, the decentralized dynamics introduce network disagreement and data heterogeneity.}
	First, network disagreement is tightly coupled with the optimization error throughout the analysis; for instance, the key one-step inequality in Lemma \ref{lem3.2} contains both a stochastic-gradient term and a consensus-error term. Second, local stochastic subgradients need not vanish near a global minimizer $\m{z}^*$. For nonsmooth objectives, the oracle may return nonzero subgradients arbitrarily close to $\m{z}^*$; even when each $f_i$ is smooth, data heterogeneity generally implies $\nabla f_i(\m{z}^*)\ne \m 0$, because the local minimizers need not coincide with a global minimizer. Consequently, the local iterates may continue to oscillate around the solution set, making their terminal values difficult to control. 

	{These observations motivate the central question of this paper: {\it can the logarithmic loss of the last iterate be removed in decentralized stochastic subgradient optimization, so that every node attains the optimal convergence rate at its last local iterate without time averaging?}}
	

	\paragraph{Contributions} 
We answer the question affirmatively by introducing a horizon-dependent multistage stepsize schedule compatible with network topology, named cooldown stepsize. We show that PDSG under this schedule attains the order-optimal $\mathcal{O}(1/\sqrt{T})$ and $\mathcal{O}(1/T)$ last-iterate rates which match the corresponding centralized lower bounds in their dependence on $T$. Specifically, we make three main contributions as follow.

First, we prescribe a unified stepsize condition for adapting to the decentralized setting, where a marked feature is that the stepsize schedule must also be compatible with network contraction. Both Condition \eqref{regular} and Assumption \ref{as_mix} provide a convenient network-independent sufficient condition for this purpose. With these conditions, our analysis clarifies how network disagreement
affects the last-iterate bound. The consensus error $\mathbb{E}\|\m{x}_i^t-\bar{\m{x}}^t\|^2 $ remains
controlled at $(\alpha^t)^2$ while its network-dependent constants
enter the coefficients lying in last-iterate bounds.

Second, as a baseline and an analytical foundation for the cooldown result, we introduce a different and unified proof compared with that of Dong et al. \cite{dong2021nonergodic} and establish last-iterate guarantees for classical PDSG under the standard stepsizes $\alpha^t=\alpha/\sqrt{t}$ and $\alpha^t=2/(\mu t)$. For every node, the expected suboptimality of the terminal local iterate is bounded by $\mathcal{O}(\log T/\sqrt{T})$ when each local objective is convex, and by $\mathcal{O}(\log T/T)$ when the aggregate objective is strongly convex. 
Notably, our analysis removes the bounded domain requirement needed by Dong et al. \cite{dong2021nonergodic} in convex cases.

Third, rather than prescribing some specific modification of standard stepsizes, we formulate a general cooldown framework in which a base stepsize sequence is combined with segmentwise scaling coefficients over a horizon-dependent partition. The construction of Jain et al.~\cite{jain2021making} is recovered as a particular instance,
whereas our conditions allow more general segment
geometries and scaling sequences. 
We know the essential requirement for stepsize construction is only
a weighted summability of the long-range stepsize ratios (Remark \ref{rem_4.4}). It is sufficient that the maximum number of
segments crossed over a lag of length $k>0$ be $o(k)$ uniformly in the horizon $T$, which intrinsically generalize the construction of Jain et al.~\cite{jain2021making} and is actually essential to the consensus analysis. By propagating objective-gap bounds between segmentwise anchors and choosing the final segment to contain only iteration $T$, we prove $\mathcal{O}(1/\sqrt{T})$ and $\mathcal{O}(1/T)$ expected optimal bounds for the last local iterate at every node.

	\paragraph{Related work}
	\noindent\textbf{Decentralized convex optimization.}
	{The convergence theory of decentralized optimization is well developed
	for smooth objectives, including EXTRA-type methods
	\cite{shi2015extra,yuan2018exact}, gradient-tracking methods
	\cite{qu2017harnessing,nedic2017achieving}, and primal--dual methods
	\cite{shi2014linear,ling2015dlm}. For nonsmooth objectives, decentralized
	subgradient methods provide a simple and widely applicable framework
	\cite{nedic2009distributed}. Their finite-time analyses, as well as
	those of many related decentralized methods \cite{koloskova2021improved,eisen2025enhancing,huang2026pareto}, commonly rely on
	time-averaged or best iterates to control oscillations due to nonsmoothness and heterogeneity. Such guarantees do not directly characterize the final local iterate available at each node.}\\
\noindent\textbf{Nonergodic decentralized methods.}
{A notable exception to averaging-based analyses is the work of Dong et
al.~\cite{dong2021nonergodic}, which study a decentralized stochastic
proximal-gradient method for nonsmooth composite optimization over
time-varying directed networks. They establish nonergodic, or
last-iterate, convergence guarantees under standard diminishing
stepsizes, while their optimal-order guarantees are obtained through
time averaging. Their framework is broader than ours in terms of
the network model and composite proximal structure, whereas our focus
is on attaining the optimal rate for the last local iterate itself in
the classical PDSG scheme. The convex guarantees most directly
comparable with ours are summarized in
Table~\ref{tab:comparison_dong}.}
\begin{table}[t]
	\centering
	\caption{Comparison of convex guarantees with
		Dong et al.~\cite{dong2021nonergodic}.}
	\label{tab:comparison_dong}
	\small
	\setlength{\tabcolsep}{4pt}
	\renewcommand{\arraystretch}{1.12}
	\begin{tabular}{lcc}
		\toprule
		Result & Rate & Bounded domain? \\
		\midrule
		
		Dong et al. \cite{dong2021nonergodic} (last iterate)
		& $\mathcal{O}(\log T/\sqrt{T})$
		& Yes
		\\
		
		Dong et al. \cite{dong2021nonergodic} (averaged iterate)
		& $\mathcal{O}(1/\sqrt{T})$
		& Yes
		\\
		
		\midrule
		
		Standard PDSG, $\alpha^t=\alpha/\sqrt{t}$
		& $\mathcal{O}(\log T/\sqrt{T})$
		& No$^\dagger$
		\\
		
		Cooldown PDSG, base $\gamma^t=\alpha/\sqrt{t}$
		& $\mathcal{O}(1/\sqrt{T})$
		& Yes
		\\
		
		Cooldown PDSG, base $\gamma^t=\alpha/\sqrt{T}$
		& $\mathcal{O}(1/\sqrt{T})$
		& No$^\dagger$
		\\
		
		\bottomrule
	\end{tabular}
	
	\vspace{1mm}
	\parbox{0.94\linewidth}{\footnotesize
		$^\dagger$For the unbounded-domain results, a finite
		initial-distance condition is imposed.}
\end{table}
{For strongly convex problems, Dong et
al.~\cite{dong2021nonergodic} state their result under strong convexity
of every local loss component, whereas our analysis requires strong
convexity only of the aggregate objective.
Nonergodic guarantees have also been established for decentralized
primal--dual and multiplier methods; see, e.g.,
Qiu and Lin~\cite{qiu2026non}, Liu et
al. \cite{liu2017projected} and Liu et al. \cite{liu2020unitary}. Such methods rely on
primal--dual dynamics and are fundamentally different from the
stochastic subgradient setting considered here.}\\
	\noindent\textbf{Centralized last-iterate convergence.}
	{For centralized stochastic subgradient methods, the gap between ergodic
	and last-iterate convergence has been studied extensively.
	Rakhlin et al.~\cite{rakhlin2012making} and Shamir and Zhang \cite{shamir2013stochastic} established last-iterate guarantees for
	standard diminishing stepsizes, while Harvey et al.
	\cite{harvey2019tight} demonstrated the tightness of the associated
	logarithmic loss. Jain et al.~\cite{jain2021making} subsequently
	introduced a multistage horizon-dependent stepsize schedule that
	recovers the optimal $\mathcal{O}(1/\sqrt{T})$ and
	$\mathcal{O}(1/T)$ rates. 
    {Defazio et al. \cite{defazio2024optimal} developed refined last-iterate analyses for horizon-dependent decay schedules, including linear decay near the terminal phase.}
     More recently, Zamani and Glineur
	\cite{zamani2025exact} derived an iteration-dependent stepsize sequence
	matching an exact worst-case lower bound and {showed that no universal stepsize sequence attains the exact optimal guarantee at every horizon.}
	Our cooldown construction adapts this last-iterate scheduling principle
	to decentralized stochastic subgradient dynamics while explicitly
	controlling post-projection network disagreement.}

	\paragraph{Organization}
	The remainder of the paper is organized as follows. Section 2 introduces the definitions, assumptions, and preliminary results. Section 3 establishes last-iterate bounds for PDSG under the standard stepsizes in the convex and strongly convex settings. Section 4 introduces the cooldown stepsize schedule and proves the corresponding optimal last-iterate rates. Section 5 presents numerical experiments illustrating the theoretical results, and Section 6 concludes the paper. 
	
	\paragraph{Notation}\label{sec:notation}
	Uppercase and lowercase boldface symbols denote matrices and vectors, respectively. For $\m v_i\in\mathbb R^p$, define $\bar{\m v}:=n^{-1}\sum_{i=1}^n\m v_i$ and $\m v:=[\m v_1;\ldots;\m v_n]\in\mathbb R^{np}$. The set $\mathcal N_i$ contains node $i$ and its neighbors. We write $\m I_p$ for the $p\times p$ identity matrix and use $\m I:=\m I_{np}$ when the dimension is clear. The symbols $\otimes$ and $(\cdot)\tr$ denote the Kronecker product and transpose. Finally, $\m M:=n^{-1}\m 1_n\m 1_n\tr\otimes\m I_p$ is the network averaging matrix.
	
	\section{Preliminaries}
	
	We encode communication over the network by a mixing matrix $\tilde{\m{W}}=[ \tilde{W}_{i j} ]\in \R^{n \times n}$, defined as follows.
	\begin{definition}\label{mix}
		(Mixing matrix $\tilde{\m{W}}$ for a given network $\mathcal{G} = \left(\mathcal{V},\mathcal{E}\right)$)
		\begin{itemize}
			\item [I.] $\tilde{\m{W}}$ is entrywise nonnegative and respects the graph topology: $\tilde{{W}}_{ij}>0$ if $j \in \mathcal{N}_i$, and $\tilde{{W}}_{ij}=0$ otherwise.
			\item [II.] $\tilde{\m{W}}$ is symmetric and doubly stochastic, i.e., $\tilde{\m{W}}=\tilde{\m{W}}\tr$ and $\tilde{\m{W}}\m{1}_n=\m{1}_n$.\\
		\end{itemize}
	\end{definition}
	Common choices for $\tilde{\m{W}}$ include the Laplacian-based constant-edge-weight matrix \cite{sayed2014diffusion} and
	the Metropolis weight matrix \cite{xiao2007distributed}.
	Let $\lambda_{i}(\tilde{\m{W}})$ denote the $i$-th largest eigenvalue of $\tilde{\m{W}}$ and define the mixing factor $\sigma:=\|\m W-\m M\|_2=\max\{|\lambda_2(\tilde{\m W})|,|\lambda_n(\tilde{\m W})|\}<1$; thus, $1-\sigma$ is the spectral gap. 
	The following properties hold.
	\begin{lemma}\label{property W}
		For $\tilde{\m{W}}$ defined in Definition~\ref{mix} and  $\m{W} :=\tilde{\m{W}}\otimes\m{I}_p$, we have
		\begin{itemize}
			\item [I.] $1=\lambda_{1}(\tilde{\m{W}})>\lambda_{2}(\tilde{\m{W}})\geq\ldots\geq\lambda_{n}(\tilde{\m{W}})>-1$;
			\item [II.] {$0\leq\|\m{W}-\m{M}\|=\sigma=\max \left\{|\lambda_{2}(\tilde{\m{W}})|, |\lambda_{n}(\tilde{\m{W}})|\right\}<1$};
			\item [III.] $\m{M}=\m{M}\m{W}=\m{W}\m{M} = \m{M}^2$;
			\item [IV.]  
			$\|\m{W}\m{x}-\m{M}\m{x}\|=\|(\m{W}-\m{M})(\m{x}-\m{M}\m{x})\| \leq \sigma\|\m{x}-\m{M}\m{x}\|$
			for any $\m{x} \in \mathbb{R}^{np}$. 
		\end{itemize}
	\end{lemma}
	
	Let $\mathcal F_t$ be the $\sigma$-algebra generated by the initial
iterates and all oracle samples observed before iteration $t$.
Hence, $\m x^t$ is $\mathcal F_t$-measurable.	
\newredfix{Throughout, we assume that the initial iterates satisfy
$\m{x}_i^1\in\Omega$ almost surely for all $i$, and that
$\mathbb{E}[\|\m{x}^1-\m{1}\otimes\m{z}^*\|^2]<\infty$.
For horizon-dependent stepsize schedules, this initial second moment
is assumed to be bounded uniformly in $T$.}
	\begin{assum}\label{as1}
		Each $f_i$ is convex on $\Omega$. There exists a constant $G > 0$ such that, for every $\m{z} \in \Omega$ and every $\hat{\m{g}}_i(\m{z}) \in \partial f_i(\m{z})$, we have $\|\hat{\m{g}}_i(\m{z})\| \leq G$. {In addition, the stochastic oracle satisfies the conditional moment bounds $\mathbb{E}[\m g_i^t\mid\mathcal F_t]=\hat{\m g}_i(\m x_i^t)\in\partial f_i(\m x_i^t)$ and $\mathbb{E}[\|\m g_i^t\|^2\mid\mathcal F_t]\le G^2$ almost surely for all $t$ and all nodes $i$.}
	\end{assum}
	
	We next record two technical lemmas used throughout the analysis.
	\begin{lemma}\label{young}
		(Young's inequality) For any two vectors $\m{v}_1,\m{v}_2\in \mathbb{R}^p$, and $\eta>0$,
		\begin{align*}
			2\m{v}_1\tr\m{v}_2 \leq \eta \|\m{v}_1\|^2+ \frac{1}{\eta}\|\m{v}_2\|^2,\quad
			\|\m{v}_1+\m{v}_2\|^2 \leq (1+\eta) \|\m{v}_1\|^2+ \left(1+\frac{1}{\eta}\right)\|\m{v}_2\|^2.
		\end{align*}
	\end{lemma}
	\begin{lemma}\label{lem:proj_variance}
		(\cite[Lemma A.2]{choi2025convergence}) Let $\Omega \subset \mathbb{R}^p$ be a closed convex set. For any set of vectors $\{\m{v}_i\}_{i=1}^n \subset \mathbb{R}^p$, let $\bar{\m{v}} = n^{-1}\sum_{i=1}^n \m{v}_i$. Then the projection operator $\mathcal{P}_\Omega$ satisfies
		$
			\sum_{i=1}^n \left\| \mathcal{P}_\Omega[\m{v}_i] - n^{-1}\sum_{k=1}^n \mathcal{P}_\Omega[\m{v}_k] \right\|^2 \leq \sum_{i=1}^n \|\m{v}_i - \bar{\m{v}}\|^2.
		$
	\end{lemma}

	\section{Last-iterate rates with standard stepsizes}
	{We first analyze PDSG under the standard diminishing stepsizes
	$\alpha^t=\alpha/\sqrt{t}$ and $\alpha^t=2/(\mu t)$. The resulting
	$\mathcal O(\log T/\sqrt T)$ and $\mathcal O(\log T/T)$
	last-local-iterate rates are consistent with the logarithmic
	nonergodic behavior known for decentralized stochastic optimization, whereas the results are established by a different and unified proof outline. 
	Here their main role is to provide a baseline for the analysis under cooldown stepsizes and to clarify the mechanisms responsible for the logarithmic loss.}

	
		{We begin with bounding the expected consensus error $\mathbb{E}[\|\m{x}_i^t - \bar{\m{x}}^t\|^2]$. 
	The ingredient is a long-range regularity condition on the stepsize sequence 
	$\{\alpha^t\}_{t \ge 1}$:
	\begin{align}\label{regular}\tag{C.1}
		(\alpha^j/\alpha^t)^2 \le C_\alpha(1+t-j)^a ~\text{for some constants}~ C_\alpha, a > 0~\text{and all}~1 \le j \le t
	\end{align}
	This condition allows the polynomial variation of the stepsizes to be absorbed by the geometric network contraction, and hence controls the entire consensus convolution (see the underlined term in \eqref{eq_unrolled_std}) without splitting the iteration horizon. It also anticipates the network-compatible scaling condition used for the cooldown schedules in Section~\ref{ord_opt}. This provides an alternative to iteration-splitting arguments such as those used in \cite{zeng2018nonconvex}.}
	
	\begin{lemma}\label{lem_std_consensus_proj}
		Suppose Assumption \ref{as1} holds and the nonincreasing stepsize sequence $\{\alpha^t\}_{t \ge 1}$ satisfies the regularity condition \eqref{regular}. Let $\{\m{x}^t\}$ be the sequence generated by PDSG. Then there exists a constant $\Cstd > 0$, independent of $t$, such that for every $t \geq 1$ and every $i \in \{1, \dots, n\}$, the expected consensus error satisfies
		$
			\mathbb{E}[\|\m{x}_i^t - \bar{\m{x}}^t\|^2] \leq \Cstd (\alpha^t)^2.
		$
	\end{lemma}
	
	\begin{proof}
		Let $Q^t := \mathbb{E}[\|\m{x}^t - \m{M}\m{x}^t\|^2] = \sum_{i=1}^n \mathbb{E}[\|\m{x}_i^t - \bar{\m{x}}^t\|^2]$. Defining the unprojected update step $\m{v}_i^t = \sum_{j=1}^n \tilde{W}_{ij}\m{x}_j^t - \alpha^t \m{g}_i^t$ and its average $\bar{\m{v}}^t = {n}^{-1}\sum_{i=1}^n \m{v}_i^t$, we have $\m{x}_i^{t+1} = \mathcal{P}_\Omega[\m{v}_i^t]$. Invoking Lemma \ref{lem:proj_variance} and using the stacked vectors $\m{x}^t$ and $\m{g}^t$ yields
		\begin{align}\label{x-mx}
			\|\m{x}^{t+1} - \m{M}\m{x}^{t+1}\|^2 &= \sum_{i=1}^n \left\| \mathcal{P}_\Omega[\m{v}_i^t] - \frac{1}{n}\sum_{k=1}^n \mathcal{P}_\Omega[\m{v}_k^t] \right\|^2 \\
			&\leq \sum_{i=1}^n \|\m{v}_i^t - \bar{\m{v}}^t\|^2 = \|(\m{W} - \m{M})\m{x}^t - \alpha^t (\m{I} - \m{M})\m{g}^t\|^2.\notag
		\end{align}
		Applying Young's inequality with a parameter $\delta > 0$, taking expectations, and invoking Lemma \ref{property W} alongside Assumption \ref{as1}, we obtain
		$
			Q^{t+1} \leq (1+\delta)\sigma^2 Q^t + \left(1+{1}/{\delta}\right)nG^2 (\alpha^t)^2.
		$
		\newredfix{If $0<\sigma<1$, choosing $\delta = {1-\sigma^2}/{(2\sigma^2)}$ yields $\lambda := (1+\delta)\sigma^2 = {(1+\sigma^2)}/{2} \in (0, 1)$ and $C_v := \left(1+{1}/{\delta}\right)nG^2 ={(1+\sigma^2)nG^2}/{(1-\sigma^2)}$. If $\sigma=0$, then $\m W=\m M$ and \eqref{x-mx} gives $Q^{t+1}\leq (\alpha^t)^2\mathbb{E}[\|(\m I-\m M)\m g^t\|^2]
			\leq nG^2(\alpha^t)^2$.
		Hence, since $Q^t\geq0$, the same recursion holds with $\lambda=1/2$ and $C_v=nG^2$. Thus, in either case, the recursion simplifies to}
		\begin{align}\label{eq_consensus_recursion}
			Q^{t+1} \leq \lambda Q^t + C_v (\alpha^t)^2.
		\end{align}
		Unrolling this recursion from $k=1$ to $t$, we obtain
		\begin{align}\label{eq_unrolled_std}
			Q^{t+1} \leq \lambda^t Q^1 + \underline{C_v \sum_{k=1}^t \lambda^{t-k} (\alpha^k)^2}.
		\end{align}
		Applying the regularity condition \eqref{regular} with $\ell := t-k$ gives
		\begin{align*}
			\sum_{k=1}^t \lambda^{t-k} (\alpha^k)^2 = (\alpha^t)^2 \sum_{k=1}^t \lambda^{t-k} \left(\frac{\alpha^k}{\alpha^t}\right)^2 \leq C_\alpha (\alpha^t)^2 \sum_{\ell=0}^{t-1} \lambda^\ell (1+\ell)^a \leq C_r (\alpha^t)^2,
		\end{align*}
		where $C_r := C_\alpha \sum_{\ell=0}^{\infty} \lambda^\ell (1+\ell)^a < \infty$ since the exponential decay of $\lambda^\ell$ asymptotically dominates the polynomial growth of $(1+\ell)^a$.
		
		Additionally, because $\lambda \in (0,1)$ and $\alpha^t$ decays at most polynomially (by \eqref{regular}), there exists $C_0 > 0$ depending on the initialization and fixed problem parameters such that $\lambda^t Q^1 \leq C_0 (\alpha^t)^2$ for all $t \geq 1$. Substituting these bounds into \eqref{eq_unrolled_std}, and noting that $\alpha^t \leq (2^a C_\alpha)^{1/2} \alpha^{t+1}$ by \eqref{regular}, we have
		$
			Q^{t+1} \leq (C_0 + C_v C_r) (\alpha^t)^2 \leq (C_0 + C_v C_r) 2^a C_\alpha (\alpha^{t+1})^2.
		$
		Since $\mathbb{E}[\|\m{x}_i^t - \bar{\m{x}}^t\|^2] \leq Q^t$ for any node $i$, setting $\mathcal{C}_{\mathrm{std}} := \max\{Q^1/(\alpha^1)^2,\, (C_0 + C_v C_r) 2^a C_\alpha\}$ ensures the stated bound holds for all $t \geq 1$.
	\end{proof}
	\begin{remark}\label{rem:regular_condition}
		Both standard stepsize schedules satisfy the regularity condition \eqref{regular}. Indeed, ${t}/{j} = 1 + {(t-j)}/{j} \le 1 + (t-j)$, and hence
		\begin{itemize}
			\item \textbf{Convex case:} The stepsize $\alpha^t = {\alpha}/{\sqrt{t}}$ yields $(\alpha^j/\alpha^t)^2 = {t}/{j} \le 1 + (t-j)$, which satisfies the condition \eqref{regular} with $C_\alpha = 1$ and $a = 1$. 
			\item \textbf{Strongly convex case:} The stepsize $\alpha^t = {2}/{(\mu t)}$ yields $(\alpha^j/\alpha^t)^2 = \left({t}/{j}\right)^2 \le (1+t-j)^2$, which satisfies the condition \eqref{regular} with $C_\alpha = 1$ and $a = 2$. 
		\end{itemize}
		Consequently, the consensus error bound derived in Lemma \ref{lem_std_consensus_proj} applies directly to both standard settings.
	\end{remark}
	\begin{remark}
		Since $\Omega$ is convex, $\bar{\m{x}}^t \in \Omega$. Assumption \ref{as1} implies that the global objective $F$ is $G$-Lipschitz continuous on $\Omega$. Thus, we have  $F(\m{x}_i^t) - F(\bar{\m{x}}^t) \leq G\|\m{x}_i^t - \bar{\m{x}}^t\|$. Taking expectations, applying the Cauchy--Schwarz inequality $\mathbb{E}[X] \leq \sqrt{\mathbb{E}[X^2]}$, and using Lemma \ref{lem_std_consensus_proj} yields 
		\begin{align}\label{eq:transfer_interface}
			\mathbb{E}[F(\m{x}_i^t)] \leq \mathbb{E}[F(\bar{\m{x}}^t)] + G\sqrt{\mathcal{C}_{\mathrm{std}}}\alpha^{t}.
		\end{align}
	\end{remark}
	
	Thus, a bound on the objective value along the network-average sequence $\{\bar{\m{x}}^t\}$ transfers to every local sequence $\{\m{x}_i^t\}$. 
	It remains to bound the objective gap at $\bar{\m{x}}^t$. 
	We begin with a one-step inequality for PDSG.
	
	\begin{lemma}\label{lem3.2}
		{	Suppose Assumption \ref{as1} holds. Let $\{\m{x}^t\}$ be the sequence generated by PDSG. For any $t_0 \ge 1$, any $\mathcal{F}_{t_0}$-measurable reference point $\m{u} \in \Omega$ \newredfix{with $\mathbb{E}[\|\m{u}-\m{z}^*\|^2]<\infty$}, and any constant $\tau > 0$, the following inequality holds for all $t \ge t_0$:
			\begin{align}\label{Fxt-Fu}
				\alpha^t \mathbb{E}[F(\bar{\m{x}}^t) - F(\m{u})] \leq & \frac{1}{2n} \mathbb{E}[\|\m{x}^{t} - \m{1}\otimes \m{u}\|^2] - \frac{1}{2n}\mathbb{E}[\|\m{x}^{t+1} - \m{1}\otimes \m{u}\|^2]  \\
				& + \frac{5\tau}{2n} \mathbb{E}[\|\m{x}^t - \m{1}\otimes\bar{\m{x}}^t\|^2] + \frac{1}{2}\left( 1 + \frac{2}{\tau} \right) G^2 (\alpha^t)^2.\notag
		\end{align}}
	\end{lemma}
	
	\begin{proof}
		Using \eqref{update_rule} and the nonexpansiveness of Euclidean projection, for the given reference point $\m{u} \in \Omega$, we have
		\begin{align*}
			\|\m{x}_i^{t+1} - \m{u}\|^2 &= \left\| \mathcal{P}_\Omega \left[ \sum_{j=1}^n \tilde{\m{W}}_{ij}\m{x}_j^t - \alpha^t \m{g}_i^t \right] - \mathcal{P}_\Omega[\m{u}] \right\|^2 
			\leq \left\| \sum_{j=1}^n \tilde{\m{W}}_{ij}\m{x}_j^t - \m{u} - \alpha^t \m{g}_i^t \right\|^2 \notag\\
			&= \left\| \sum_{j=1}^n \tilde{\m{W}}_{ij}(\m{x}_j^t - \m{u}) \right\|^2 - 2\alpha^t \left\langle \sum_{j=1}^n \tilde{\m{W}}_{ij}\m{x}_j^t - \m{u}, \m{g}_i^t \right\rangle + (\alpha^t)^2 \|\m{g}_i^t\|^2.
		\end{align*}
		Summing this over $i=1, \ldots, n$ and invoking Lemma \ref{property W}, we obtain
		\begin{align}\label{eq_sum_x_update}
			\sum_{i=1}^n \|\m{x}_i^{t+1} - \m{u}\|^2 \leq & \sum_{i=1}^n \|\m{x}_i^t - \m{u}\|^2 - 2\alpha^t \sum_{i=1}^n \langle \m{x}_i^t - \m{u}, \m{g}_i^t \rangle  \\
			& - 2\alpha^t \sum_{i=1}^n \left\langle \sum_{j=1}^n \tilde{\m{W}}_{ij}\m{x}_j^t - \m{x}_i^t, \m{g}_i^t \right\rangle + (\alpha^t)^2 \sum_{i=1}^n \|\m{g}_i^t\|^2.\notag
		\end{align}
		\redfix{For $t \ge t_0$, both $\m{x}_i^t$ and the $\mathcal{F}_{t_0}$-measurable reference point $\m{u}$ are $\mathcal{F}_t$-measurable. The tower property therefore gives
			\begin{align*}
				\mathbb{E}[\langle \m{x}_i^t - \m{u}, \m{g}_i^t \rangle] &= \mathbb{E} \big[ \mathbb{E}[\langle \m{x}_i^t - \m{u}, \m{g}_i^t \rangle \mid \mathcal{F}_t] \big] 
				= \mathbb{E} \big[ \langle \m{x}_i^t - \m{u}, \mathbb{E}[\m{g}_i^t \mid \mathcal{F}_t] \rangle \big] \\
				&= \mathbb{E}[\langle \m{x}_i^t - \m{u}, \hat{\m{g}}_i(\m{x}_i^t) \rangle] \geq \mathbb{E}[f_i(\m{x}_i^t) - f_i(\m{u})].
		\end{align*}}
		Applying Young's inequality with an arbitrary constant $\tau > 0$ gives
		\begin{align*}
			-2\alpha^t \sum_{i=1}^n\left\langle \sum_{j=1}^n \tilde{\m{W}}_{ij}\m{x}_j^t - \m{x}_i^t, \m{g}_i^t \right\rangle 
			\leq &\tau \|(\m{W} - \m{I})\m{x}^t\|^2 + \frac{(\alpha^t)^2}{\tau} \|\m{g}^t\|^2\\
			\stackrel{\text{Lemma }\ref{property W}}{\leq} &4\tau \|\m{x}^t - \m{1}\otimes\bar{\m{x}}^t\|^2+\frac{(\alpha^t)^2}{\tau} \|\m{g}^t\|^2.
		\end{align*}
		Substituting the preceding two inequalities into \eqref{eq_sum_x_update} and dividing by $2n$ gives
		\begin{align}\label{eq_intermediate_gap}
			&\alpha^t \mathbb{E}\left[ \frac{1}{n}\sum_{i=1}^n \big(f_i(\m{x}_i^t) - f_i(\m{u})\big) \right] \leq  \frac{1}{2n} \mathbb{E}[\|\m{x}^{t}-\m{1}\otimes \m{u}\|^2] \\
			&- \frac{1}{2n}\mathbb{E}[\|\m{x}^{t+1}-\m{1}\otimes \m{u}\|^2]  
			+ \frac{2\tau}{n} \mathbb{E}[\|\m{x}^t-\m{1}\otimes\bar{\m{x}}^t\|^2] + \frac{1}{2}\left( 1 + \frac{1}{\tau} \right) G^2 ({\alpha^t})^2.\notag
		\end{align}
		Finally, by the convexity of each $f_i$, we have $f_i(\m{x}_i^t) \geq f_i(\bar{\m{x}}^t) + \langle \hat{\m{g}}_i(\bar{\m{x}}^t), \m{x}_i^t - \bar{\m{x}}^t \rangle$. Summing over $i$ and then taking expectations, we obtain
		\begin{align}\label{eq_intermediate_gap1}
			\mathbb{E}\left[ \frac{1}{n}\sum_{i=1}^n f_i(\m{x}_i^t) \right] \geq \mathbb{E}[F(\bar{\m{x}}^t)] - \frac{1}{n} \sum_{i=1}^n \mathbb{E}\big[ \| \hat{\m{g}}_i(\bar{\m{x}}^t) \| \|\m{x}_i^t - \bar{\m{x}}^t\| \big].
		\end{align}
		Applying Young's inequality again, we obtain
		\begin{align}\label{eq_intermediate_gap2}
			\alpha^t \mathbb{E}\big[ \| \hat{\m{g}}_i(\bar{\m{x}}^t) \| \|\m{x}_i^t - \bar{\m{x}}^t\| \big] \leq \frac{\tau}{2} \mathbb{E}[\|\m{x}_i^t - \bar{\m{x}}^t\|^2] + \frac{(\alpha^t)^2}{2\tau} \mathbb{E}[\|\hat{\m{g}}_i(\bar{\m{x}}^t)\|^2].
		\end{align}
		Substituting \eqref{eq_intermediate_gap1} and \eqref{eq_intermediate_gap2} into \eqref{eq_intermediate_gap} completes the proof.
	\end{proof}
	
	Combining the one-step inequality with Lemma \ref{lem_std_consensus_proj}, \eqref{Fxt-Fu} can be written as
	\begin{align}\label{Fxt-Fu-compact}
		&	\alpha^t \mathbb{E}[F(\bar{\m{x}}^t)-F(\m{u})] \\
		\leq &\frac{1}{2n} (\mathbb{E}[\|\m{x}^{t}-\m{1}\otimes \m{u}\|^2] - \mathbb{E}[\|\m{x}^{t+1}-\m{1}\otimes \m{u}\|^2] )+ K_{std}({\alpha^t})^2,\notag
	\end{align}
	where $K_{std} := \big({5\tau} \Cstd + \left(1+{2}/{\tau}\right)G^2\big)/2$.
	Setting $\m{u}=\m{z}^*$ in \eqref{Fxt-Fu-compact}, we get
	\begin{align}\label{Fxt-Fz}
		&	\alpha^t(\mathbb{E}[F(\bar{\m{x}}^t)-F^*]) \\
		\leq &\frac{1}{2n} (\mathbb{E}[\|\m{x}^{t}-\m{1}\otimes\m{z}^*\|^2]-\mathbb{E}[\|\m{x}^{t+1}-\m{1}\otimes\m{z}^*\|^2] )+ K_{std}({\alpha^t})^2.\notag
	\end{align}
	Summing the preceding inequalities over $t=1,\ldots,T$ gives
	\begin{align}\label{sumFxt-Fz}
		\sum_{t=1}^{T} \alpha^t\left(\mathbb{E}[F(\bar{\m{x}}^t)]-F^*\right) \leq \frac{1}{2n}\mathbb{E}\|\m{x}^{1}-\m{1}\otimes\m{z}^*\|^2 + \sum_{t=1}^{T} K_{std}({\alpha^t})^2.
	\end{align}
	
	While \eqref{sumFxt-Fz} controls a weighted cumulative objective gap along the $\{\bar{\m{x}}^t\}$, last-iterate analysis in nonsmooth optimization requires a reduction from this aggregate bound to the terminal objective gap. The next lemma provides such a bridge.
	
	\begin{lemma}\label{lem3.5}
		(\cite[Lemma 4]{defazio2024optimal}) Let $\{\nu^t\}$ be a positive, nonincreasing sequence, and let $q^t\geq0$ for every $t$. Then
		\begin{align}\label{lem3.5_eq1}
			\nu^Tq^T \leq \frac{1}{T} \sum_{t=1}^{T}\nu^tq^t+\sum_{t=1}^{T-1}\frac{1}{(t+1)t} \sum_{k=T-t+1}^{T} \nu^k (q^k-q^{(T-t)}).
		\end{align}
	\end{lemma}
	We are now ready to state the convex last-iterate guarantee for PDSG. 	
	
	\subsection{Convex objectives} 
	\begin{theorem}\label{thm_last_iterate_cvx}
		Suppose Assumption \ref{as1} holds. Let $T\ge2$, choose $\alpha^t = \alpha/\sqrt{t}$ for some constant $\alpha > 0$, and let $\{\m{x}^t\}$ be generated by PDSG. Then, for every node $i \in \{1, \dots, n\}$, the expected suboptimality of the last local iterate $\m{x}_i^T$ satisfies
		$
			\mathbb{E}[F(\m{x}_i^T)] - F^* \leq \mathcal{O}\left({\log T}/{\sqrt{T}}\right).
		$
		More precisely, the following explicit bound holds:
		\begin{align*}
			\mathbb{E}[F(\m{x}_i^T)] - F^* \leq \frac{\frac{1}{n}\mathbb{E}\|\m{x}^{1}-\m{1}\otimes\m{z}^*\|^2}{2\alpha \sqrt{T}} + \frac{(3K_{std} + \Cstd)\alpha(1+\ln T)}{\sqrt{T}} + \frac{G\sqrt{\Cstd}\alpha}{\sqrt{T}}.
		\end{align*}
	\end{theorem}
	
	\begin{proof}
		Applying Lemma \ref{lem3.5} with $q^t=\mathbb{E}[F(\bar{\m{x}}^t)]-F^*$ and $\nu^t=\alpha^t$, we obtain
		\begin{align}\label{eq1+}
			&\alpha^T\left(\mathbb{E}[F(\bar{\m{x}}^T)]-F^*\right) \leq \frac{1}{T}\sum_{t=1}^{T} \alpha^t\left(\mathbb{E}[F(\bar{\m{x}}^t)]-F^*\right)\\
			&+ \sum_{t=1}^{T-1}\frac{1}{(t+1)t} \sum_{k=T-t+1}^{T} \alpha^k \mathbb{E}[F(\bar{\m{x}}^k)-F(\bar{\m{x}}^{(T-t)})].\notag
		\end{align}
		Note that in \eqref{eq1+}, we can extend the starting index from $k=T-t+1$ to $k=T-t$ without changing the summation, since the term evaluates to zero when $k=T-t$. To bound the inner sum, fix the outer index $t$ and condition on $\mathcal{F}_{T-t}$. We apply the one-step descent inequality \eqref{Fxt-Fu-compact} from step $k=T-t$ to $T$ using the $\mathcal{F}_{T-t}$-measurable reference point $\m{u}=\bar{\m{x}}^{T-t}$. Summing these inequalities, taking expectations, and using the tower property yields
		\begin{align}\label{eq1++}
			&\sum_{k=T-t}^T\alpha^k \mathbb{E}[F(\bar{\m{x}}^k)-F(\bar{\m{x}}^{T-t})] \\
			\leq &\frac{1}{2n} \mathbb{E}[\|\m{x}^{T-t}-\m{1}\otimes\bar{\m{x}}^{T-t}\|^2]+\sum_{k=T-t}^{T} K_{std}({\alpha^k})^2.\notag
		\end{align}
		By Lemma \ref{lem_std_consensus_proj}, we have $(2n)^{-1}\mathbb{E}[\|\m{x}^{T-t}-\m{1}\otimes\bar{\m{x}}^{T-t}\|^2] \leq {\Cstd} (\alpha^{T-t})^2/{2}$. Substituting this bound into \eqref{eq1++} and then into \eqref{eq1+} yields
		\begin{align}\label{eq2_new}
			&\alpha^T\left(\mathbb{E}[F(\bar{\m{x}}^T)]-F^*\right) \leq \frac{1}{T}\sum_{t=1}^{T} \alpha^t\left(\mathbb{E}[F(\bar{\m{x}}^t)]-F^*\right)\\
			&+\sum_{t=1}^{T-1}\frac{1}{(t+1)t} \left( \frac{\Cstd}{2} (\alpha^{T-t})^2 + \sum_{k=T-t}^{T} K_{std}({\alpha^k})^2 \right).\notag
		\end{align}
		Substituting $\alpha^t=\alpha/\sqrt{t}$ into \eqref{sumFxt-Fz} gives
		\begin{align}\label{sum_bound1_new}
			\sum_{t=1}^{T} \alpha^t\left(\mathbb{E}[F(\bar{\m{x}}^t)]-F^*\right) &\leq \frac{1}{2n}\mathbb{E}\|\m{x}^{1}-\m{1}\otimes\m{z}^*\|^2+ \alpha^2 K_{std}\sum_{t=1}^{T}\frac{1}{t}\\
			&\leq \frac{1}{2n}\mathbb{E}\|\m{x}^{1}-\m{1}\otimes\m{z}^*\|^2 + K_{std}\alpha^2(1+\ln T).\notag
		\end{align}
		Observe that $\sum_{k=T-t+1}^{T}\frac{ 1}{k} \leq \int_{T-t}^T \frac{1}{x} \mbox{d}x=\ln \frac{T}{T-t}=\ln \left(1+\frac{t}{T-t}\right) \leq \frac{t}{T-t}$. Hence,
		\begin{align*}
			\sum_{k=T-t}^{T}\frac{ 1}{k} = \frac{1}{T-t}+\sum_{k=T-t+1}^{T}\frac{ 1}{k} \leq \frac{1}{T-t} + \frac{t}{T-t} = \frac{t+1}{T-t}.
		\end{align*}
		Since ${(t+1)}/{(T-t)} > {1}/{(T-t)}$, we can bound the double summation in \eqref{eq2_new} by
		\begin{align}\label{sum_bound2_new}
			&\sum_{t=1}^{T-1}\frac{1}{(t+1)t} \left( \frac{\Cstd \alpha^2}{2(T-t)} + K_{std}\alpha^2\sum_{k=T-t}^{T}\frac{1}{k} \right) \\
			\leq & \alpha^2 \left(K_{std} + \frac{\Cstd}{2}\right) \sum_{t=1}^{T-1}\frac{1}{(t+1)t} \frac{t+1}{T-t} = \alpha^2 \left(K_{std} + \frac{\Cstd}{2}\right) \sum_{t=1}^{T-1}\frac{1}{t(T-t)} \notag\\
			= &\frac{\alpha^2}{T} \left(K_{std} + \frac{\Cstd}{2}\right) \sum_{t=1}^{T-1}\left(\frac{1}{t}+\frac{1}{T-t}\right) \leq \frac{2\alpha^2 (K_{std} + \Cstd/2)(1+\ln T)}{T}.\notag
		\end{align}
		Dividing \eqref{eq2_new} by $\alpha^T = \alpha/\sqrt{T}$ and using \eqref{sum_bound1_new} and \eqref{sum_bound2_new} yields the following bound for the network average at time $T$:
		\begin{align}\label{xbar-z}
			\mathbb{E}[F(\bar{\m{x}}^T)]-F^* \leq  \frac{\frac{1}{n}\mathbb{E}\|\m{x}^{1}-\m{1}\otimes\m{z}^*\|^2}{2\alpha \sqrt{T}} + \frac{(3K_{std} + \Cstd)\alpha(1+\ln T)}{\sqrt{T}}.
		\end{align}
		Finally, combining \eqref{xbar-z} with the transfer bound \eqref{eq:transfer_interface} gives the stated last-iterate convergence rate for the convex case.
	\end{proof}
	
	\subsection{Strongly convex objectives}
	\begin{lemma}
		For any $k\in\{1,\ldots,T-1\}$, any $\mathcal{F}_{T-k}$-measurable reference point $\m{u} \in \Omega$ \newredfix{with $\mathbb{E}[\|\m{u}-\m{z}^*\|^2]<\infty$},
		\begin{align}\label{lem4.1 eq1}
			&\sum_{t=T-k}^T \mathbb{E}[F(\bar{\m{x}}^t)-F(\m{u})]\leq \frac{\mu(T-k)}{4n}\mathbb{E}[\|\m{x}^{T-k}-\m{1}\otimes\m{u}\|^2]\\
			&+\frac{\mu}{4n}\sum_{t=T-k+1}^T \mathbb{E}[\|\m{x}^{t}-\m{1}\otimes\m{u}\|^2]+\frac{2K_{std}}{\mu}\sum_{t=T-k}^T\frac{1}{t}.\notag
		\end{align}
	\end{lemma}
	\begin{proof}
		Dividing \eqref{Fxt-Fu-compact} by $\alpha^t$ and summing over $t=T-k, \ldots, T$ yields
		\begin{align*}
			&\sum_{t=T-k}^T \mathbb{E}[F(\bar{\m{x}}^t)-F(\m{u})]\leq \frac{1}{2n \alpha^{T-k}}\mathbb{E}[\|\m{x}^{T-k}-\m{1}\otimes\m{u}\|^2]\\
			&+\frac{1}{2n}\sum_{t=T-k+1}^T \mathbb{E}[\|\m{x}^{t}-\m{1}\otimes\m{u}\|^2]\left(\frac{1}{\alpha^t}-\frac{1}{\alpha^{t-1}}\right)+\sum_{t=T-k}^T K_{std} \alpha^t.
		\end{align*}
		Substituting $\alpha^t={2}/{(\mu t)}$ gives $\left({1}/{\alpha^t}-{1}/{\alpha^{t-1}}\right) = {\mu}/{2}$, which proves the claim.
	\end{proof}
	
	\begin{lemma}\label{lem_sc_dist}
		Suppose that the global objective $F$ is $\mu$-strongly convex \newredfix{on $\Omega$} and that $\alpha^t = {2}/{(\mu t)}$. Let $D^t := {(n)^{-1}}\mathbb{E}[\|\m{x}^{t}-\m{1}\otimes\m{z}^*\|^2]$ denote the average expected squared distance to the optimum, and define the constant {$Q_{sc} := \max(2G^2+\Cstd, 2K_{std}+2\Cstd)$}. Then, for every $t \geq 1$, the following bound holds:
		\begin{align}\label{sc-conv}
			D^t \leq {4 Q_{sc}}/{(\mu^2 t)}.
		\end{align}
	\end{lemma}
	\begin{proof}
		By the strong convexity of $F$ and the fact that $\m{z}^*$ minimizes $F$ \newredfix{over the convex set $\Omega$}, we obtain $F(\bar{\m{x}}^t)-F^* \geq {\mu}/{2} \|\bar{\m{x}}^t-\m{z}^*\|^2$. \newredfix{Indeed, for any $\theta\in(0,1)$, strong convexity and the optimality of $\m{z}^*$ give
		\[
		F^* \leq F((1-\theta)\m{z}^*+\theta\bar{\m{x}}^t)
		\leq (1-\theta)F^*+\theta F(\bar{\m{x}}^t)-\frac{\mu}{2}\theta(1-\theta)\|\bar{\m{x}}^t-\m{z}^*\|^2.
		\]
		Letting $\theta\downarrow0$ yields the stated inequality.} Note that the average distance decomposes as $D^t = \mathbb{E}[\|\bar{\m{x}}^t-\m{z}^*\|^2] + {n}^{-1}\mathbb{E}[\|\m{x}^t - \m{1}\otimes\bar{\m{x}}^t\|^2]$. Using the consensus error bound from Lemma \ref{lem_std_consensus_proj}, we have
		\newredfix{\begin{align*}
		\mathbb{E}[\|\bar{\m{x}}^t-\m{z}^*\|^2]
		= D^t - \frac{1}{n}\mathbb{E}[\|\m{x}^t - \m{1}\otimes\bar{\m{x}}^t\|^2]
		\geq D^t - \Cstd(\alpha^t)^2.
		\end{align*}}
		Combining \eqref{Fxt-Fz} with these bounds gives
		\begin{align*}
			&D^{t+1} \leq D^t - 2\alpha^t \mathbb{E}[F(\bar{\m{x}}^t)-F^*] + 2K_{std}(\alpha^t)^2 \\
			&\leq D^t - \mu\alpha^t (D^t - \Cstd(\alpha^t)^2) + 2K_{std}(\alpha^t)^2 
			\leq (1-\mu\alpha^t)D^t + (2K_{std} + \mu \Cstd \alpha^t) (\alpha^t)^2.
		\end{align*}
		Since $\alpha^t = {2}/{(\mu t)}$, we have $\mu\alpha^t \leq 2$ for all $t \geq 1$. Thus, $2K_{std} + \mu \Cstd \alpha^t \leq 2K_{std} + 2\Cstd \leq Q_{sc}$. The recursion simplifies to
		\begin{align}\label{eq_Dt_recursion}
			D^{t+1} \leq \left(1-{2}/{t}\right)D^t + 4Q_{sc}/{(\mu^2 t^2)}.
		\end{align}
		
		To establish \eqref{sc-conv}, first consider \(t=1,2\). Define
		$
		\hat{\m h}(\bar{\m x}^t)
		:=\frac1n\sum_{i=1}^n\hat{\m g}_i(\bar{\m x}^t)
		\in\partial F(\bar{\m x}^t).
		$
		Strong convexity and Cauchy--Schwarz give
		\(
		\|\hat{\m h}(\bar{\m x}^t)\|
		\ge(\mu/2)\|\bar{\m x}^t-\m z^*\|
		\).
		Assumption~\ref{as1} also gives
		\(
		\|\hat{\m h}(\bar{\m x}^t)\|
		\le n^{-1}\sum_i\|\hat{\m g}_i(\bar{\m x}^t)\|\le G
		\).
		Consequently, for every \(t\ge1\),
		\begin{equation}\label{u_b}
			\mathbb{E}[\|\bar{\m{x}}^t-\m{z}^*\|^2]
			\le \frac{4}{\mu^2}
			\mathbb{E}[\|\hat{\m h}(\bar{\m x}^t)\|^2]
			\le \frac{4G^2}{\mu^2}.
		\end{equation}
		
		For $t=1$, combining \eqref{u_b} with Lemma \ref{lem_std_consensus_proj} and using $\alpha^1 = {2}/{\mu}$ gives
		\begin{align*}
			D^1
			&=\mathbb{E}[\|\bar{\m{x}}^1-\m{z}^*\|^2]
			+\frac1n\mathbb{E}[\|\m{x}^1-\m{1}\otimes\bar{\m{x}}^1\|^2]\\
			&\le \frac{4G^2}{\mu^2}+\Cstd(\alpha^1)^2
			=\frac{4(G^2+\Cstd)}{\mu^2}
			\le\frac{4Q_{sc}}{\mu^2}.
		\end{align*}
		which establishes \eqref{sc-conv} for $t=1$. 
		Similarly, for $t=2$, with $\alpha^2 = {1}/{\mu}$, we have
		\begin{align*}
			D^2 \leq \frac{4G^2}{\mu^2} + \Cstd(\alpha^2)^2 = \frac{4G^2+\Cstd}{\mu^2} \leq \frac{4 Q_{sc}}{\mu^2 \cdot 2},
		\end{align*}
		Thus, \eqref{sc-conv} also holds for $t=2$.
		
		We now proceed by induction for $t\ge3$. Evaluating the recursion \eqref{eq_Dt_recursion} at $t=2$, we obtain
		\begin{align*}
			D^3 \leq \left(1-\frac{2}{2}\right)D^2 + Q_{sc} \frac{4}{\mu^2 (2)^2} = 0 + \frac{Q_{sc}}{\mu^2} = \frac{3 Q_{sc}}{3\mu^2} \leq \frac{4 Q_{sc}}{\mu^2 \cdot 3}.
		\end{align*}
		This establishes the base case $t=3$. Now assume that \eqref{sc-conv} holds for some integer $k \geq 3$. For $t=k+1$, we deduce
		\begin{align*}
			D^{k+1} \leq \left(1-\frac{2}{k}\right)\frac{4 Q_{sc}}{\mu^2 k}+ \frac{4 Q_{sc}}{\mu^2 k^2} 
			= \frac{4 Q_{sc}}{\mu^2}\left(\frac{k-1}{k^2}\right) \leq \frac{4 Q_{sc}}{\mu^2 (k+1)}.
		\end{align*}
		This completes the induction.
	\end{proof}
	
	The next two lemmas bound, respectively, the two terms on the right-hand side of \eqref{lem3.5_eq1} when it is applied with the constant weight sequence $\nu^t\equiv1$.

	\begin{lemma}\label{lem_3.7}
		For any $k\in\{1,\ldots,T-1\}$,
		\begin{align}\label{eq5}
			\sum_{t=T-k}^T \mathbb{E}[F(\bar{\m{x}}^t)-F(\bar{\m{x}}^{T-k})]\leq \frac{\Cstd+2K_{std}+(4Q_{sc}+2K_{std})k}{\mu(T-k)}.
		\end{align}
	\end{lemma}
	\begin{proof}
		Setting $\m{u}=\bar{\m{x}}^{T-k}$ in \eqref{lem4.1 eq1}, we bound the terms on the right-hand side separately.
		Invoking Lemma \ref{lem_std_consensus_proj} together with  $\alpha^{T-k} = {2}/{\mu(T-k)}$, we have
		\begin{align}\label{eq_residual_boundary}
			{\frac{\mu(T-k)}{4n} \mathbb{E}[\|\m{x}^{T-k} - \m{1}\otimes\bar{\m{x}}^{T-k}\|^2] }
			\le \frac{\mu(T-k)}{4} \Cstd (\alpha^{T-k})^2 
			= \frac{\Cstd}{\mu(T-k)}.
		\end{align}
		For any $t \geq T-k$, we have
		\begin{align}\label{lem_3.7_eq1}
			\frac{1}{n}\mathbb{E}[\|\m{x}^{t}-\m{1}\otimes\bar{\m{x}}^{T-k}\|^2] &\leq \frac{2}{n}\mathbb{E}[\|\m{x}^{t}-\m{1}\otimes\m{z}^*\|^2] + 2\mathbb{E}[\|\bar{\m{x}}^{T-k}-\m{z}^*\|^2] \\
			&\leq 2 D^t + 2 D^{T-k} \stackrel{\text{Lemma }\ref{lem_sc_dist}}{\leq} \frac{8 Q_{sc}}{\mu^2} \left(\frac{1}{t}+\frac{1}{T-k}\right).\notag		\end{align}
		Summing \eqref{lem_3.7_eq1} over $t=T-k+1,\ldots,T$ and multiplying by $\mu/4$ gives
		\begin{align*}
			{\frac{\mu}{4}\sum_{t=T-k+1}^T \frac{1}{n}\mathbb{E}[\|\m{x}^{t}-\m{1}\otimes\bar{\m{x}}^{T-k}\|^2]} 
			\leq \frac{2Q_{sc}}{\mu}\left(\sum_{t=T-k+1}^T \frac{1}{t} + \frac{k}{T-k}\right)\leq \frac{4 Q_{sc} k}{\mu(T-k)}.
		\end{align*}	
		Substituting this inquality and \eqref{eq_residual_boundary}  into \eqref{lem4.1 eq1} and using $\sum_{t=T-k}^{T}\frac{1}{t} \leq \frac{k+1}{T-k}$ yields the stated inequality.
	\end{proof}

	\begin{lemma}\label{lem_sc_average_bound}
		Suppose that the assumptions and notation of Lemma \ref{lem_sc_dist} hold. Then, for all $T \ge 1$,
		\begin{align}\label{eq_sc_average_bound}
			\frac{1}{T}\sum_{t=1}^{T} \big(\mathbb{E}[F(\bar{\m{x}}^t)]-F^*\big) \leq \frac{(Q_{sc}+2K_{std})(1+\ln T)}{\mu T}.
		\end{align}
	\end{lemma}
	\begin{proof}
		Dividing \eqref{Fxt-Fz} by $\alpha^t$ gives for all $t \ge 1$,
		\begin{align*}
			\mathbb{E}[F(\bar{\m{x}}^t)-F^*] \leq \frac{1}{2\alpha^t} D^t - \frac{1}{2\alpha^t} D^{t+1} + K_{std} \alpha^t.
		\end{align*}
		Summing over $t=1,\ldots,T$ and writing $A_t := {1}/{(2\alpha^t)}$, we obtain
		$
			\sum_{t=1}^{T}\big(\mathbb{E}[F(\bar{\m{x}}^t)]-F^*\big) \leq \sum_{t=1}^{T}\big(A_tD^t-A_tD^{t+1}\big) + K_{std}\sum_{t=1}^{T}\alpha^t.
		$
		Rearranging the telescoping sum gives
		\begin{align*}
			&\sum_{t=1}^{T}\big(A_tD^t-A_tD^{t+1}\big) = A_1D^1 - A_TD^{T+1} \\
			&+ \sum_{t=2}^{T}(A_t-A_{t-1})D^t \leq A_1D^1 + \sum_{t=2}^{T}(A_t-A_{t-1})D^t,
		\end{align*}
		Since $\alpha^t=\frac{2}{\mu t}$, we have $A_1=\frac{1}{2\alpha^1}=\frac{\mu}{4}$ and $A_t-A_{t-1}=\frac{1}{2}\left(\frac{1}{\alpha^t}-\frac{1}{\alpha^{t-1}}\right)=\frac{\mu}{4}$ for all $t\ge2$. Hence
		$	A_1D^1 + \sum_{t=2}^{T}(A_t-A_{t-1})D^t = \frac{\mu}{4}\sum_{t=1}^{T}D^t.$
		Applying Lemma \ref{lem_sc_dist}, which gives $D^t \leq {4Q_{sc}}/{\mu^2 t}$ for \emph{all} $t\ge1$, we get
		\begin{align*}
			\frac{\mu}{4}\sum_{t=1}^{T}D^t \leq \frac{Q_{sc}}{\mu}\sum_{t=1}^{T}\frac{1}{t} \leq \frac{Q_{sc}}{\mu}(1+\ln T).
		\end{align*}
		Similarly, $K_{std}\sum_{t=1}^{T}\alpha^t = ({2K_{std}}/{\mu})\sum_{t=1}^{T}{1}/{t} \leq ({2K_{std}}/{\mu})(1+\ln T)$. Combining the two bounds and dividing by $T$ yields \eqref{eq_sc_average_bound}.
	\end{proof}

	\begin{theorem}\label{thm_last_iterate_sc}
		Suppose Assumption \ref{as1} holds and the global objective $F$ is $\mu$-strongly convex \newredfix{on $\Omega$}. Let $T\ge2$, choose $\alpha^t = {2}/{\mu t}$, and let $\{\m{x}^t\}$ be generated by PDSG. Then, for every node $i \in \{1, \dots, n\}$, the expected suboptimality of 
		$\m{x}_i^T$ satisfies
		\begin{align*}
			\mathbb{E}[F(\m{x}_i^T)] - F^* \leq \mathcal{O}\left(\frac{\log T}{T}\right).
		\end{align*}
		More precisely, the following explicit bound holds:
		\begin{align*}
			\mathbb{E}[F(\m{x}_i^T)] - F^* \leq \frac{1}{\mu T}\Big[(9Q_{sc}+10K_{std}+2\Cstd)(1+\ln T) + 2G\sqrt{\Cstd}\Big],
		\end{align*}
	\end{theorem}
	\begin{proof}
		Apply Lemma \ref{lem3.5} with the constant weight sequence $\nu^t\equiv1$ and $q^t:=\mathbb{E}[F(\bar{\m{x}}^t)]-F^*\ge0$. This gives
		\begin{align}\label{eq_sc_lem35_applied}
			&\mathbb{E}[F(\bar{\m{x}}^T)]-F^* \leq \underbrace{\frac{1}{T}\sum_{t=1}^{T} \big(\mathbb{E}[F(\bar{\m{x}}^t)]-F^*\big)}_{\text{(I)}}\\ 
			+& \underbrace{\sum_{j=1}^{T-1}\frac{1}{j(j+1)} \sum_{k=T-j}^{T} \big(\mathbb{E}[F(\bar{\m{x}}^k)]-\mathbb{E}[F(\bar{\m{x}}^{T-j})]\big)}_{\text{(II)}}.\notag
		\end{align}
		Term (I) is bounded by Lemma \ref{lem_sc_average_bound}; it remains to bound term (II).
		\begin{align*}
			\text{Term (II)} \stackrel{\eqref{eq5}}{\leq}  \frac{1}{\mu}\sum_{j=1}^{T-1} \frac{\Cstd+2K_{std}+(4Q_{sc}+2K_{std})j}{j(j+1)(T-j)}.
		\end{align*}
		Using $1/(j+1)\le1$ and the harmonic bound $\sum_{j=1}^{T-1}1/j\le1+\ln T$, we obtain
		\begin{align*}
			\sum_{j=1}^{T-1} \frac{1}{j(j+1)(T-j)} \leq \sum_{j=1}^{T-1} \frac{1}{j(T-j)} = \frac{1}{T}\sum_{j=1}^{T-1}\left(\frac{1}{j}+\frac{1}{T-j}\right) \leq \frac{2(1+\ln T)}{T},
		\end{align*}
		\begin{align*}
			\sum_{j=1}^{T-1}\frac{1}{(j+1)(T-j)} = \frac{1}{T+1}\sum_{j=1}^{T-1}\left(\frac{1}{j+1}+\frac{1}{T-j}\right) \leq \frac{2(1+\ln T)}{T+1} \leq \frac{2(1+\ln T)}{T}.
		\end{align*}
		Combining the two estimates gives
		\begin{align}\label{eq_sc_termII}
			\text{Term (II)} \leq \frac{2(\Cstd+4K_{std}+4Q_{sc})(1+\ln T)}{\mu T}.
		\end{align}
		Substituting \eqref{eq_sc_average_bound} and \eqref{eq_sc_termII} into \eqref{eq_sc_lem35_applied}, we obtain
		\begin{align}\label{eq_sc_final_avg}
			\mathbb{E}[F(\bar{\m{x}}^T)]-F^* \leq \frac{(9Q_{sc}+10K_{std}+2\Cstd)(1+\ln T)}{\mu T}.
		\end{align}
		Combining \eqref{eq_sc_final_avg} with \eqref{eq:transfer_interface} and $\alpha^T={2}/{(\mu T)}$ yields the stated $\mathcal{O}\left({\log T}/{T}\right)$ rate.
	\end{proof}
	
	\section{A general cooldown framework for order-optimal last-iterate rates}\label{ord_opt}
	
	{We now develop a general class of horizon-dependent multistage stepsize schedules, consisting of three
components: a base sequence, a horizon-dependent partition, and a
sequence of segmentwise scaling coefficients. This formulation allows
us to identify separately the conditions required by the optimization
dynamics (Assumption \ref{as3}) and those required by network mixing (Assumption \ref{as_mix}).}

	Let
	$0=T_0<T_1<\cdots<T_K<T_{K+1}=T$ and
	$\Delta_i:=T_{i+1}-T_i,
	$
	where \(K=K(T)\) is the number of stepsize reductions, and let
	\(\{p_i\}_{i=0}^K\subset(0,1]\) be the corresponding scaling coefficients.
	
	\begin{assum}\label{as3}
		\textbf{(Cooldown structure)}
		There exist constants \(0<c_-\le c_+<1\), \(D\ge1\),
		\(\rho\in(0,1)\), and \(M>0\), independent of \(T\), such that
		\begin{align*}
			&1=p_0>p_1>\cdots>p_K>0,\qquad
			p_i\ge \rho p_{i-1}\ (1\le i\le K),\qquad
			\sum_{i=0}^K p_i\le M,\\
			&c_-T\le \Delta_0=T_1\le c_+T,\qquad
			\Delta_K=1,\qquad
			\Delta_i\le D\Delta_{i+1},
			\qquad 0\le i\le K-1.
		\end{align*}
	\end{assum}
	
	\begin{assum}\label{as_mix}
		\textbf{(Network-compatible scaling)}
		Let \(i(t)=i\) when \(T_i<t\le T_{i+1}\). There exist constants
		\(C_p\ge1\) and \(s\ge0\), independent of \(T\), such that
		$
			{p_{i(j)}}/{p_{i(t)}}\le C_p(1+t-j)^s,
			\, 1\le j\le t\le T$.
	\end{assum}
	
	Assumption~\ref{as_mix} controls long-range stepsize ratios, which are
	convolved with the geometric network contraction in the consensus
	analysis. It does not prescribe the number of cooldown segments. The
	following proposition records two useful sufficient constructions.
	
	\begin{proposition}\label{prop_admissible_schedules}
		\textnormal{(i)} Suppose Assumption~\ref{as3} holds and, in addition,
		\(\Delta_i\ge d\Delta_{i+1}\) for \(0\le i\le K-2\), where \(d>1\).
		Then Assumption~\ref{as_mix} holds with
		$
		C_p=\rho^{-2},\,
		s={\ln(1/\rho)}/{\ln d}.
		$
		\textnormal{(ii)} For \(T\ge4\), define
		\[
		K=\lceil\log_2T\rceil,\qquad
		T_i=T-\lceil T2^{-i}\rceil\ (0\le i\le K),\qquad
		T_{K+1}=T,\qquad p_i=2^{-i}.
		\]
		The dyadic-ceiling schedule of Jain et al.~\cite{jain2021making} satisfies Assumptions~\ref{as3}
		and~\ref{as_mix} with
		\[
		c_-=2/5,\quad c_+=1/2,\quad D=4,\quad
		\rho=1/2,\quad M=2,\quad C_p=4,\quad s=1.
		\]
	\end{proposition}
	
	\begin{proof}
		For part \textnormal{(i)}, let \(m=i(t)-i(j)\). If \(m\ge2\), the
		geometric lower bound on the segment lengths gives
		\(t-j>d^{m-2}\). Hence
		$
		{p_{i(j)}}/{p_{i(t)}}\le \rho^{-m}
		\le \rho^{-2}(t-j)^{\ln(1/\rho)/\ln d};
		$
		the cases \(m\le1\) are immediate.
		
		For part \textnormal{(ii)}, \(\Delta_0=\lfloor T/2\rfloor\), and the
		ceiling identity \(2\lceil a\rceil-1\le\lceil2a\rceil\le2\lceil
		a\rceil\) gives \(\Delta_i\le4\Delta_{i+1}\); the remaining conditions
		in Assumption~\ref{as3} follow directly. To verify
		Assumption~\ref{as_mix}, take \(j<t\), set
		\(m=i(t)-i(j)\), and write \(k=t-j\). The claim is immediate for
		\(m\le1\). If \(m\ge2\), set \(x=T2^{-i(t)}>1/2\). Then
		$
		k>T_{i(t)}-T_{i(j)+1}
		=\lceil2^{m-1}x\rceil-\lceil x\rceil
		\ge2^{m-2}.
		$
		Therefore \(p_{i(j)}/p_{i(t)}=2^m\le4k\le4(1+k)\).
	\end{proof}

    \begin{remark}\label{rem_4.4}
{Assumption~\ref{as_mix} is a network-independent sufficient condition
for Lemma~\ref{lem_consensus}. The proof in fact only requires
 $
    S:=\sup_{T,t}
    \sum_{k=0}^{t-1}
    \lambda^k
    \left(
        {\alpha^{t-k}}/{\alpha^t}
    \right)^2
    <\infty, \, 0<\lambda<1.
$
To interpret this condition, define $ N(k):=
    \sup_{T}\sup_{k<t\le T}
    \big(i(t)-i(t-k)\big)$,
the maximum number of segments that can be crossed over a
lag of length $k$. Under Assumption~\ref{as3}, ${p_{i(t-k)}}/{p_{i(t)}}
    \le \rho^{-N(k)}$,
and, for $\gamma^t=ct^{-q}$, ${\gamma^{t-k}}/{\gamma^t}
    \le (1+k)^q$.
Hence,
$\left(
        {\alpha^{t-k}}/{\alpha^t}
    \right)^2 \leq  \rho^{-2N(k)}(1+k)^{2q}$.
{Therefore, if $N(k)=o(k)$ uniformly in $T$, then for every $\widetilde{\lambda}\in(\lambda,1)$, there exists $k_0$ such that
$
\lambda^k\rho^{-2N(k)}(1+k)^{2q}
\le
\widetilde{\lambda}^{\,k},
\, k\ge k_0.
$
Hence $S$ is finite. We use the polynomial condition in Assumption 4.2 because it is independent of the network contraction factor and directly covers the dyadic-ceiling construction \cite{jain2021making}.}
}
\end{remark}

	The cooldown stepsize schedule is defined by
	\begin{equation}\label{mod_step}
		\alpha^t=p_i\gamma^t,
		\qquad T_i<t\le T_{i+1},\quad 0\le i\le K,
	\end{equation}
	where
	\[
	\gamma^t=ct^{-q},\qquad
	(c,q)=
	\begin{cases}
		(\alpha,\frac12), & \text{for convex objectives},\\
		(\frac{2}{\mu},1), & \text{for strongly convex objectives}.
	\end{cases}
	\]
	The sequence \(\{\gamma^t\}\) is nonincreasing and, for every
	\(\epsilon\in(0,1)\) and \(\epsilon T\le t\le T\),
	$
	\gamma^T
	=\left({t}/{T}\right)^q\gamma^t
	\ge \epsilon^q\gamma^t
	\ge \epsilon\gamma^t.
	$
	\begin{remark}
		The initial segment retains the standard stepsize scale, while the later
		segments reduce the update magnitude. The bound
		\(\sum_i p_i\le M\) controls the segmentwise optimization errors, and
		Assumption~\ref{as_mix} controls their interaction with network mixing.
		The singleton final segment makes its anchor the prescribed output at
		time \(T\). No bound on \(K\) is required.
	\end{remark}

	The proof telescopes over segmentwise anchors
	\(\tau_0,\ldots,\tau_K\), with \(\tau_K=T\). The initial phase
	provides an ergodic bound at \(\tau_0\), and Lemma~\ref{lem_descent}
	bounds the \(i\)th anchor increment by
	\(\mathcal{O}(p_i\gamma^T)\). Thus \(\sum_i p_i\le M\), rather than
	the number of segments, controls the accumulated cooldown error.

	\begin{lemma}\label{lem_consensus}
		Suppose Assumptions \ref{as1}, \ref{as3}, and \ref{as_mix} hold. Let $\{\m{x}^t\}$ be the sequence generated by PDSG using the cooldown stepsize schedule defined in \eqref{mod_step}. Then there exists a constant $\Cmod > 0$ such that, for every $T \geq 2$, every $1 \leq t \leq T$, and every node $i \in \{1, \dots, n\}$, the expected consensus error satisfies
		\begin{align}\label{lem_consensus_eq1}
			\mathbb{E}[\|\m{x}_i^t - \bar{\m{x}}^t\|^2] \leq \Cmod (\alpha^t)^2.
		\end{align}
	\end{lemma}
	\begin{proof}
		See Appendix \ref{proof_lem_consensus}.
	\end{proof}

	{\begin{lemma}\label{lem_gap_sum}
			Suppose Assumptions \ref{as1}, \ref{as3}, and \ref{as_mix} hold. Let $\{\m{x}^t\}$ be generated by PDSG using the cooldown stepsize defined in \eqref{mod_step}. For any indices $1 \leq t_0 \leq t_1 \leq T$ and any $\mathcal{F}_{t_0}$-measurable reference point $\m{u} \in \Omega$ \newredfix{with $\mathbb{E}[\|\m{u}-\m{z}^*\|^2]<\infty$}, the following inequality holds:
			\begin{align}\label{lem6.2_eq_u}
				\sum_{t=t_0}^{t_1} 2\alpha^t (\mathbb{E}[F(\bar{\m{x}}^t) - F(\m{u})]) \leq \frac{1}{n}\mathbb{E}[\|\m{x}^{t_0}-\m{1}\otimes\m{u}\|^2] + 2K_{mod} \sum_{t=t_0}^{t_1} (\alpha^t)^2,
			\end{align}
			where $K_{mod} = ({5\tau}\Cmod + \left(1+{2}/{\tau}\right)G^2)/2$.
			In particular, setting $\m{u} = \bar{\m{x}}^{t_0}$, we have
			\begin{align}\label{lem6.2_eq}
				\sum_{t=t_0}^{t_1} 2\alpha^t (\mathbb{E}[F(\bar{\m{x}}^t) - F(\bar{\m{x}}^{t_0})]) \leq Q_{mod} \sum_{t=t_0}^{t_1} (\alpha^t)^2,
			\end{align}
			where $Q_{mod} = 2K_{mod} + \Cmod$.
		\end{lemma}
		
		\begin{proof}
			As in Lemma \ref{lem3.2}, the one-step descent inequality under the cooldown schedule holds for every $\mathcal{F}_{t_0}$-measurable reference point $\m{u} \in \Omega$ and $t \ge t_0$. Applying Lemma \ref{lem_consensus} gives
			\begin{align}\label{Fxt-Fu-mod}
				&2\alpha^t(\mathbb{E}[F(\bar{\m{x}}^t)-F(\m{u})])\\ \leq & \frac{1}{n}\mathbb{E}[\|\m{x}^{t}-\m{1}\otimes\m{u}\|^2]-\frac{1}{n}\mathbb{E}[\|\m{x}^{t+1}-\m{1}\otimes\m{u}\|^2]+2K_{mod}({\alpha^t})^2.\notag
			\end{align}
			Summing \eqref{Fxt-Fu-mod} from $t = t_0$ to $t_1$ and dropping the nonpositive term yields \eqref{lem6.2_eq_u}.
			
			Because the network average $\bar{\m{x}}^{t_0}$ is $\mathcal{F}_{t_0}$-measurable, setting $\m{u} = \bar{\m{x}}^{t_0}$ in \eqref{lem6.2_eq_u} gives
			\begin{align*}
				\sum_{t=t_0}^{t_1} 2\alpha^t (\mathbb{E}[F(\bar{\m{x}}^t)-F(\bar{\m{x}}^{t_0})]) \leq \frac{1}{n}\mathbb{E}[\|\m{x}^{t_0}-\m{1}\otimes\bar{\m{x}}^{t_0}\|^2] + 2K_{mod} \sum_{t=t_0}^{t_1} (\alpha^t)^2.
			\end{align*}
			Applying \eqref{lem_consensus_eq1} yields \eqref{lem6.2_eq}.	
	\end{proof}}

	For each segment, define an index attaining the minimum expected objective value:
	\begin{align*}
		\tau_0 = \arg\min_{\lfloor T_1/2 \rfloor < t \leq T_1} \mathbb{E}[F(\bar{\m{x}}^t)],\quad\tau_i = \arg\min_{T_i < t \leq T_{i+1}} \mathbb{E}[F(\bar{\m{x}}^t)],~i \in \{1,\ldots,K\},
	\end{align*}
	These deterministic indices serve as anchors for comparing successive cooldown segments. The following lemma bounds the objective-value increment between consecutive anchors.
	
	\begin{lemma}\label{lem_descent}
		Suppose Assumptions \ref{as1}, \ref{as3}, and \ref{as_mix} hold. Let $\{\m{x}^t\}$ be the sequence generated by PDSG using the cooldown stepsize schedule defined in \eqref{mod_step}. Then, for every $i \in \{1,\ldots,K-1\}$, 
		\begin{align*}
			\mathbb{E}[F(\bar{\m{x}}^{\tau_{i+1}}) - F(\bar{\m{x}}^{\tau_i})] \leq \frac{(D+1) Q_{mod}}{2 \rho} \kappa_{\gamma} p_{i} \gamma^{T},
		\end{align*}
		where $\kappa_{\gamma} = {1}/{c_-}$ for convex objectives, and $\kappa_{\gamma} = {1}/{c_-^2}$ for strongly convex objectives.
	\end{lemma}
	\begin{proof}
		Define $\Delta^+_i := \mathbb{E}[F(\bar{\m{x}}^{\tau_{i+1}})-F(\bar{\m{x}}^{\tau_{i}})]$. If $\Delta^+_i \leq 0$, the claim is immediate; hence assume $\Delta^+_i > 0$. Setting $t_0=\tau_i$ and $t_1=T_{i+2}$ in Lemma \ref{lem_gap_sum} yields
		\begin{align}\label{proof_eq_base}
			\sum_{t=\tau_i}^{T_{i+2}} 2\alpha^t (\mathbb{E}[F(\bar{\m{x}}^t)-F(\bar{\m{x}}^{\tau_i})]) \leq Q_{mod} \sum_{t=\tau_i}^{T_{i+2}} (\alpha^t)^2.
		\end{align}
		Since $\alpha^t \leq \alpha^{T_i+1} = p_i \gamma^{T_i+1}$ for all $t \in [\tau_i, T_{i+2}]$, and the definition of $\tau_i$ implies $\tau_i \geq T_i + 1$, we have $T_{i+2}-\tau_i+1 \leq T_{i+2}-T_i$ and thus the RHS of \eqref{proof_eq_base} satisfies
		\begin{align}\label{proof_rhs_bound}
			\text{RHS} \leq Q_{mod} (T_{i+2}-T_i) (p_i \gamma^{T_i+1})^2.
		\end{align}
		For the LHS of \eqref{proof_eq_base}, restricting the summation domain to $t \in [T_{i+1}+1, T_{i+2}]$ and using  $\mathbb{E}[F(\bar{\m{x}}^t)] \geq \mathbb{E}[F(\bar{\m{x}}^{\tau_{i+1}})]$ and $\alpha^t \geq p_{i+1} \gamma^{T_{i+2}}$ yields
		\begin{align}\label{proof_lhs_bound}
			\text{LHS} \geq \sum_{t=T_{i+1}+1}^{T_{i+2}} 2 p_{i+1} \gamma^{T_{i+2}} \Delta^+_i = 2 (T_{i+2}-T_{i+1}) p_{i+1} \gamma^{T_{i+2}} \Delta^+_i.
		\end{align}
		Combining \eqref{proof_rhs_bound} and \eqref{proof_lhs_bound} with \eqref{proof_eq_base} gives
	    $
			2 (T_{i+2}-T_{i+1}) p_{i+1} \gamma^{T_{i+2}} \Delta^+_i \leq Q_{mod} (T_{i+2}-T_i) p_i^2 (\gamma^{T_i+1})^2,
		$
		which implies
		\begin{align}\label{proof_eq_iso}
			\Delta^+_i \leq \frac{Q_{mod}}{2} \left( \frac{T_{i+2}-T_i}{T_{i+2}-T_{i+1}} \right) \left( \frac{p_i^2}{p_{i+1}} \right) \Lambda(T),
		\end{align}
		where $\Lambda(T) = \frac{(\gamma^{T_i+1})^2}{\gamma^{T_{i+2}}}$. By Assumption \ref{as3}, we have $\frac{T_{i+2}-T_i}{T_{i+2}-T_{i+1}} = 1 + \frac{T_{i+1}-T_i}{T_{i+2}-T_{i+1}} \leq 1+D$ and $\frac{p_i^2}{p_{i+1}} \leq \frac{p_i}{\rho}$. 
		
		We now evaluate $\Lambda(T)$ in the two cases.
		
		\noindent$\bullet$ \textbf{Convex case:} $\gamma^t = {\alpha}/{\sqrt{t}} \implies \Lambda(T) = \alpha {\sqrt{T_{i+2}}}/{T_i+1}$. Since $i \geq 1$, using $T_{i+2} \leq T$ and $T_i+1 > T_i \geq c_- T$ yields
		\begin{align*}
			\Lambda(T) \leq \alpha \frac{\sqrt{T}}{c_- T} = \frac{1}{c_-} \left( \frac{\alpha}{\sqrt{T}} \right) = \frac{1}{c_-} \gamma^T.
		\end{align*}
		
		\noindent$\bullet$ \textbf{Strongly convex case:}  $\gamma^t = {2}/{\mu t} \implies \Lambda(T) = {2T_{i+2}}/{(\mu(T_i+1)^2)}$. Since $i \geq 1$, we have $T_i \geq T_1 \geq c_- T$. Combining this with $T_{i+2} \leq T$ yields
		\begin{align*}
			\Lambda(T) \leq \frac{2}{\mu} \frac{T}{(T_i+1)^2} \leq \frac{2}{\mu} \frac{T}{(c_- T)^2} = \frac{1}{c_-^2} \left(\frac{2}{\mu T}\right) = \frac{1}{c_-^2} \gamma^T.
		\end{align*}
		Substituting these bounds into \eqref{proof_eq_iso} completes the proof.
	\end{proof}
	\begin{corollary}\label{cor_first_jump}
		Under the assumptions of Lemma~\ref{lem_descent},
		\begin{align}\label{cor_first_jump_eq}
			\mathbb{E}[F(\bar{\m{x}}^{\tau_1})-F(\bar{\m{x}}^{\tau_0})]
			\le \frac{(D+1)Q_{mod}}{2\rho}\kappa_0\gamma^T,
		\end{align}
		where \(\kappa_0=2\kappa_\gamma\) in the convex case and
		\(\kappa_0=4\kappa_\gamma\) in the strongly convex case.
	\end{corollary}
	\begin{proof}
		Repeat the proof of Lemma~\ref{lem_descent} with
		\(t_0=\tau_0\) and \(t_1=T_2\). Since
		\(\tau_0>\lfloor T_1/2\rfloor\),
		$
		{(T_2-\tau_0+1)}/{\Delta_1}
		\le 1+{(\lceil T_1/2\rceil)}/{\Delta_1}
		\le 1+D.
		$
		Moreover,
		\((\gamma^{\tau_0})^2/\gamma^{T_2}\le\kappa_0\gamma^T\),
		by \(T_1\ge c_-T\). Substitution in
		\eqref{proof_eq_iso} proves the claim; the same argument covers
		\(K=1\).
	\end{proof}
	

	\begin{theorem}\label{thm_last_iterate}
		Suppose Assumptions \ref{as1}, \ref{as3}, and \ref{as_mix} hold. Let $T\ge2$ and $\{\m{x}^t\}$ be the sequence generated by PDSG using the cooldown stepsize schedule defined in \eqref{mod_step}. Then, for any node $i \in \{1, \dots, n\}$, the following bounds hold for the expected suboptimality of the last local iterate $\m{x}_i^T$:
		
		\vspace{1mm}\noindent\textbf{Convex case:} If each local objective $f_i$ is convex, $\operatorname{diam}(\Omega) \leq R$ for some $R>0$, and the base stepsize is $\gamma^t = {\alpha}/{\sqrt{t}}$ for some $\alpha > 0$, then
		$\mathbb{E}[F(\m{x}_i^T)] - F^* \leq {C_{cvx}}/{\sqrt{T}}$,
		where the constant $C_{cvx}$ is defined as
		\begin{align*}
			C_{cvx} := \frac{R^2}{\alpha \sqrt{c_-}} + \alpha \left[ \frac{2K_{mod} \ln 3}{\sqrt{c_-}} + \frac{(D+1)Q_{mod}}{2} \left( \frac{2}{c_- \rho} + \frac{M}{\rho c_-} \right) \right] + \alpha G\sqrt{\mathcal{C}_{\mathrm{mod}}}.
		\end{align*}
		
		\vspace{1mm}\noindent\textbf{Strongly convex case:} If the global objective $F$ is $\mu$-strongly convex \newredfix{on $\Omega$} and the base stepsize is $\gamma^t = {2}/{(\mu t)}$, then
		$\mathbb{E}[F(\m{x}_i^T)] - F^* \leq {C_{sc}}/{T}$,
		where $Q_{sc}^{mod} := \max(2G^2+\Cmod, 2K_{mod}+2\Cmod)$ and the constant $C_{sc}$ is defined as
		\begin{align*}
			C_{sc} := \newredfix{\frac{2 Q_{sc}^{mod}(1+\ln 2) + 4K_{mod}}{\mu c_-}} + \frac{2}{\mu} \left[ \frac{(D+1)Q_{mod}}{2} \left( \frac{4}{c_-^2 \rho} + \frac{M}{\rho c_-^2} \right) + G\sqrt{\mathcal{C}_{\mathrm{mod}}} \right].
		\end{align*}
	\end{theorem}
	
	\begin{proof}
		By the $G$-Lipschitz continuity of $F$ on $\Omega$, Lemma \ref{lem_consensus}, and $\alpha^T = p_K \gamma^T \leq \gamma^T$, we have
		\begin{align}\label{proof_step1}
			\mathbb{E}[F(\m{x}_i^T)] - \mathbb{E}[F(\bar{\m{x}}^T)] \leq G \sqrt{\mathbb{E}[\|\m{x}_i^T - \bar{\m{x}}^T\|^2]} \leq G\sqrt{\Cmod}\alpha^T \leq G\sqrt{\Cmod}\gamma^T.
		\end{align}
		Using $\tau_K = T$, the telescoping decomposition over the indices $\{\tau_{i}\}$ yields
		\begin{align}\label{proof_step2}
			\mathbb{E}[F(\bar{\m{x}}^T)] = &\mathbb{E}[F(\bar{\m{x}}^{\tau_0})] + \left( \mathbb{E}[F(\bar{\m{x}}^{\tau_1})] - \mathbb{E}[F(\bar{\m{x}}^{\tau_0})] \right) \\
			&+ \sum_{j=1}^{K-1} \left( \mathbb{E}[F(\bar{\m{x}}^{\tau_{j+1}})] - \mathbb{E}[F(\bar{\m{x}}^{\tau_j})] \right).\notag
		\end{align}
		Applying Lemma \ref{lem_descent} and using $\sum_{j=1}^{K-1} p_j \leq M$ from Assumption \ref{as3} give
		\begin{align}\label{proof_j_ge1}
			&\sum_{j=1}^{K-1} \left( \mathbb{E}[F(\bar{\m{x}}^{\tau_{j+1}})] - \mathbb{E}[F(\bar{\m{x}}^{\tau_j})] \right)\\
			\leq &\frac{(D+1)Q_{mod}\kappa_\gamma}{2\rho} \gamma^T \sum_{j=1}^{K-1} p_j 
			\leq \frac{M(D+1)Q_{mod}\kappa_\gamma}{2\rho} \gamma^T.\notag
		\end{align}
		Combining \eqref{cor_first_jump_eq} and \eqref{proof_j_ge1} in \eqref{proof_step2} yields
		\begin{align}\label{proof_step3}
			\mathbb{E}[F(\bar{\m{x}}^T)] - \mathbb{E}[F(\bar{\m{x}}^{\tau_0})] \leq \frac{(D+1)Q_{mod}}{2} \left( \frac{\kappa_0}{\rho} + \frac{M\kappa_\gamma}{\rho} \right) \gamma^T.
		\end{align}
		It remains to bound $\mathbb{E}[F(\bar{\m{x}}^{\tau_0})] - F^*$ using the suffix of the initial phase.
		
		{	\noindent$\bullet$ \textbf{Convex case:} $\gamma^t = {\alpha}/{\sqrt{t}}$, $\kappa_0 = {2}/{c_-}$, and $\kappa_\gamma = {1}/{c_-}$. Note that
			\begin{align*}
				\mathbb{E}[F(\bar{\m{x}}^{\tau_0})]-F^*
				=\min_{\lfloor T_1/2\rfloor<t\le T_1}
				\mathbb{E}[F(\bar{\m{x}}^t)-F^*]
				\le\frac{1}{T_1-\lfloor T_1/2\rfloor}
				\sum_{t=\lfloor T_1/2\rfloor+1}^{T_1}
				\mathbb{E}[F(\bar{\m{x}}^t)-F^*].
			\end{align*}
			We bound the initial-phase error over the suffix window $t \in [\lfloor T_1/2 \rfloor + 1, T_1]$. Summing the relation \eqref{Fxt-Fu-mod} over $[\lfloor T_1/2 \rfloor + 1, T_1]$ yields
			\begin{align*}
				\sum_{t=\lfloor T_1/2 \rfloor + 1}^{T_1} 2\alpha^t (\mathbb{E}[F(\bar{\m{x}}^t)] - F^*) \leq \frac{1}{n}\mathbb{E}[\|\m{x}^{\lfloor T_1/2 \rfloor + 1} - \m{1}\otimes\m{z}^*\|^2] + 2K_{mod} \sum_{t=\lfloor T_1/2 \rfloor + 1}^{T_1} (\alpha^t)^2.
			\end{align*}
			For the left-hand side, because $\alpha^t$ is nonincreasing, $\alpha^t \geq {\alpha}/{\sqrt{T_1}}$, and hence
			\begin{align*}
				\text{LHS} \geq 2 \frac{\alpha}{\sqrt{T_1}} (T_1 - \lfloor T_1/2 \rfloor) \mathbb{E}[F(\bar{\m{x}}^{\tau_0})- F^*] \geq \alpha \sqrt{T_1} \mathbb{E}[F(\bar{\m{x}}^{\tau_0})- F^*].
			\end{align*}
			For the right-hand side, we have
			\begin{align*}
				\sum_{t=\lfloor T_1/2 \rfloor + 1}^{T_1} (\alpha^t)^2 \leq \alpha^2 \int_{\lfloor T_1/2 \rfloor}^{T_1} \frac{1}{x} \text{d}x = \alpha^2 \ln \left( \frac{T_1}{\lfloor T_1/2 \rfloor} \right) \leq \alpha^2 \ln 3.
			\end{align*}
			Since $\operatorname{diam}(\Omega) \leq R$, $\frac{1}{n}\mathbb{E}[\|\m{x}^{\lfloor T_1/2 \rfloor + 1} - \m{1}\otimes\m{z}^*\|^2] \leq  R^2$. Combining these two bounds yields
			$
				\mathbb{E}[F(\bar{\m{x}}^{\tau_0})] - F^* \leq \left( R^2 + 2K_{mod}\alpha^2\ln 3 \right)/{\alpha \sqrt{T_1}} .
			$
			Since $T_1 \ge c_- T$, we obtain
			\begin{align*}
				\mathbb{E}[F(\bar{\m{x}}^{\tau_0})] - F^* \leq \left( \frac{R^2}{\alpha^2 \sqrt{c_-}} + \frac{2K_{mod} \ln 3}{\sqrt{c_-}} \right) \frac{\alpha}{\sqrt{T}} = \left( \frac{R^2}{\alpha^2 \sqrt{c_-}} + \frac{2K_{mod} \ln 3}{\sqrt{c_-}} \right) \gamma^T,
			\end{align*}
			Together with \eqref{proof_step3} and \eqref{proof_step1}, this yields the stated inequality in the convex case.
			
			\noindent$\bullet$ \textbf{Strongly convex case:} $\gamma^t = {2}/{(\mu t)}$, $\kappa_0 = {4}/{c_-^2}$, and $\kappa_\gamma = {1}/{c_-^2}$. Lemma \ref{lem_base_sc_optimal} bounds the minimum expected objective value over the initial suffix:
			\begin{align*}
				\mathbb{E}[F(\bar{\m{x}}^{\tau_0})] - F^* \newredfix{\leq \frac{2 Q_{sc}^{mod}(1+\ln 2) + 4K_{mod}}{\mu T_1} 
				\leq \frac{2 Q_{sc}^{mod}(1+\ln 2) + 4K_{mod}}{2 c_-} \gamma^T},
			\end{align*}
			where the last inequality follows from $T_1 \ge c_- T$ and ${1}/{T} = {\mu \gamma^T}/{2}$. Combining this bound with \eqref{proof_step3} and \eqref{proof_step1} gives the stated result in the strongly convex case.}
	\end{proof}
	
	For convex objectives, the bounded-domain assumption can be replaced by a finite initial-distance condition when the horizon-dependent constant base stepsize $\gamma^t \equiv \gamma := {\alpha}/{\sqrt{T}}$ is used.
	\begin{corollary}\label{cor_unbounded_domain}
		Suppose Assumptions \ref{as1}, \ref{as3}, and \ref{as_mix} hold,
		\(T\ge2\), and each \(f_i\) is convex. Let
		\(D_0:=n^{-1}\mathbb{E}[\|\m{x}^1-\m{1}\otimes\m{z}^*\|^2]\).
		If \(\gamma^t\equiv\alpha/\sqrt{T}\) in \eqref{mod_step}, then
		there exists a constant \(\widetilde C_{cvx}>0\), independent of
		\(T\), such that, for every node \(i\),
		$
			\mathbb{E}[F(\m{x}_i^T)]-F^*
			\le {\widetilde C_{cvx}}/{\sqrt T}.
		$
	\end{corollary}
	
	\begin{proof}
		Let $Q^t := \mathbb{E}[\|\m{x}^t-\m{M}\m{x}^t\|^2]$. When $\gamma^t\equiv\gamma$, Assumption~\ref{as_mix} gives
		$
		\left({\alpha^j}/{\alpha^t}\right)^2
		\le C_p^2(1+t-j)^{2s}.
		$
		Unrolling the recursion
		$Q^{t+1}\le\lambda Q^t+\bar C(\alpha^t)^2$, given in the proof of Lemma \ref{lem_consensus}, therefore yields 
		\begin{align}\label{eq_exact_consensus}
			Q^t \le \lambda^{t-1}Q^1 + \Chor(\alpha^t)^2, \quad \forall t \ge 1,
		\end{align}
		where $\Gamma_s(\lambda):=\sum_{k=0}^{\infty}(1+k)^{2s}\lambda^k<\infty$ and $\Chor := \bar{C}C_p^2\Gamma_s(\lambda)$. Because the arithmetic mean $\bar{\m{x}}^1$ minimizes $\sum_i\|\m{x}_i^1-\m{c}\|^2$ over $\m{c}$, we have $Q^1 \le nD_0$. 
		
		For $t \ge \lfloor T_1/2 \rfloor + 1 \ge c_-T/2$, Assumption~\ref{as_mix} with $j=1$ gives $p_{i(t)}\ge C_p^{-1}t^{-s}\ge C_p^{-1}T^{-s}$. Hence $\alpha^t\ge \alpha/(C_pT^{s+1/2})$ and $\lambda^{t-1}\le A_\lambda T^{-2s-1}\le (A_\lambda C_p^2/\alpha^2)(\alpha^t)^2$, where $A_\lambda:=\sup_{T\ge2}T^{2s+1}\lambda^{c_-T/2-1}<\infty$. Substituting this bound into \eqref{eq_exact_consensus} yields 
		\begin{align}\label{eq_tail_consensus}
			Q^t \le \mathcal{C}_{\text{hor}}^+(\alpha^t)^2, \quad \forall t \ge \lfloor T_1/2 \rfloor + 1,
		\end{align}
		where $\mathcal{C}_{\text{hor}}^+ := \Chor + nD_0A_\lambda C_p^2/\alpha^2$.
        
		 Summing \eqref{Fxt-Fu} over  $[1, T_1]$ with $\m{u}=\m{z}^*$ with $\alpha^t \equiv \gamma$ yields
		\begin{align}\label{proof_cor_sum1}
			2\gamma \sum_{t=1}^{T_1} (\mathbb{E}[F(\bar{\m{x}}^t)] - F^*) \le D_0 + \frac{5\tau}{n}\sum_{t=1}^{T_1}Q^t + \Big(1+\frac{2}{\tau}\Big)G^2 T_1\gamma^2.
		\end{align}
		By \eqref{eq_exact_consensus}, $\sum_{t=1}^{T_1}Q^t \le \frac{nD_0}{1-\lambda} + \Chor T_1\gamma^2$. Substituting this into \eqref{proof_cor_sum1} and dividing by $2\gamma T_1 \ge 2\gamma c_- T$ gives
        $
        2\gamma\sum_{t=1}^{T_1}
    \big(\mathbb{E}[F(\bar{\m{x}}^t)]-F^*\big)
    \leq
    \widehat D_0
    +2K_{\rm hor}T_1\gamma^2,
       $
		where $\hat{D}_0 := D_0\big(1+\frac{5\tau}{1-\lambda}\big)$ and $K_{\text{hor}} := \frac{5\tau}{2}\Chor + \frac{1}{2}\big(1+\frac{2}{\tau}\big)G^2$. By the definition of $\tau_0$,

        \begin{align*}
    \mathbb{E}[F(\bar{\m{x}}^{\tau_0})]-F^*
\leq
    \frac{1}{L_0}
    \sum_{t=\lfloor T_1/2\rfloor+1}^{T_1}
    \big(\mathbb{E}[F(\bar{\m{x}}^t)]-F^*\big)
\leq
    \frac{1}{L_0}
    \sum_{t=1}^{T_1}
    \big(\mathbb{E}[F(\bar{\m{x}}^t)]-F^*\big),
\end{align*}
where $ L_0:=T_1-\lfloor T_1/2\rfloor\ge \frac{T_1}{2}$. Thus,
		\begin{align}\label{proof_cor_base}
			\mathbb{E}[F(\bar{\m{x}}^{\tau_0})] - F^* \le \frac{\hat{D}_0}{\alpha c_-\sqrt{T}} + \frac{2K_{\text{hor}}\alpha}{\sqrt{T}},
		\end{align}
        
		For $0 \le i \le K-1$, since $\gamma^t \equiv \gamma$, we have $\Lambda_i(T) = \gamma^2/\gamma = \gamma$ for all $i$. Furthermore, $\tau_i \ge \lfloor T_1/2 \rfloor + 1$ ensures that \eqref{eq_tail_consensus} holds. Thus, Lemma \ref{lem_descent} applies for all $i \ge 0$, yielding
$\mathbb{E}[F(\bar{\m{x}}^{\tau_{i+1}}) - F(\bar{\m{x}}^{\tau_i})] \le \frac{(D+1)Q_{mod}^+}{2\rho}\, p_i \gamma$,
		where $Q_{mod}^+ := 2K_{\text{hor}}^+ + \mathcal{C}_{\text{hor}}^+$ and $K_{\text{hor}}^+ := ({5\tau}\mathcal{C}_{\text{hor}}^+ + \big(1+{2}/{\tau}\big)G^2)/2$. Summing over $i$ and invoking $\sum_{i=0}^{K-1} p_i \le M$ gives
$\mathbb{E}[F(\bar{\m{x}}^T)] - \mathbb{E}[F(\bar{\m{x}}^{\tau_0})] \le \frac{(D+1)M}{2\rho}Q_{mod}^+\,\gamma$.
		By $G$-Lipschitz continuity of $F$ and $\alpha^T \le \gamma$, \eqref{eq_tail_consensus} implies
		\begin{align}\label{proof_cor_transfer}
			\mathbb{E}[F(\m{x}_i^T)] - \mathbb{E}[F(\bar{\m{x}}^T)] \le G\sqrt{\mathcal{C}_{\text{hor}}^+}\,\alpha^T \le G\sqrt{\mathcal{C}_{\text{hor}}^+}\,\gamma.
		\end{align}
		Combining \eqref{proof_cor_base}--\eqref{proof_cor_transfer} and substituting $\gamma = {\alpha}/{\sqrt{T}}$ yields the desired result.
	\end{proof}

	\section{Numerical experiments}\label{sec:numerics}
	
	We use two mechanism checks. The first isolates the logarithmic last-iterate penalty and its removal by cooldown; the second compares last-iterate accuracy while monitoring network disagreement.
	
	Unless stated otherwise, we use $n=8$ nodes on a ring with $\tilde W_{ii}=\tilde W_{i,i\pm1}=1/3$ (hence $\sigma\approx0.805$), initialize all iterates at the origin, and let every node take one stochastic-gradient sample and one communication step per iteration. 
    For concreteness, we instantiate our general cooldown framework using
the construction of
Jain et al.~\cite{jain2021making}, that is
	\[
	K=\lceil\log_2T\rceil,\qquad
	T_\ell=T-\lceil T2^{-\ell}\rceil\ (0\le\ell\le K),\qquad
	T_{K+1}=T,
	\]
	with $\alpha^t=p_\ell\gamma^t$ for $T_\ell<t\le T_{\ell+1}$ and $p_\ell=2^{-\ell}$. Proposition~\ref{prop_admissible_schedules} shows that this construction satisfies Assumptions~\ref{as3} and~\ref{as_mix}, including when consecutive terminal segments have the same length. The base sequence $\gamma^t$ is specified below. When an exact optimum is unavailable, $F^*$ in the plots denotes the smallest objective value obtained in a longer centralized reference run.
	
\begin{figure}[htbp]
		\centering
		\includegraphics[width=.8\linewidth]{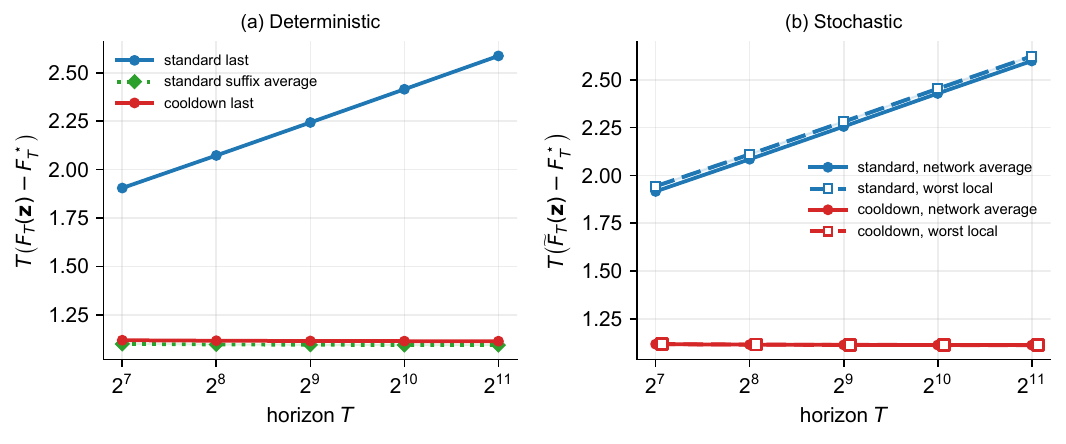}
		\caption{Lower-bound instance of Harvey et al.~\cite{harvey2019tight}. Deterministic case: the three curves are $T(F_T(\bar{\mathbf x}^{T})-F_T^*)$ for the standard and cooldown last iterates, and $T(F_T(\widehat{\mathbf x}_{1/2}^{T})-F_T^*)$ for the standard half-suffix average. Stochastic case: the four curves are the network-average and worst-local last-iterate errors under the standard and cooldown schedules, namely $T(\widetilde F_T(\bar{\mathbf x}^{T})-F_T^*)$ and $T(\widetilde F_T(\mathbf x_{i^\star}^{T})-F_T^*)$ for each schedule, where $i^\star\in\arg\max_i\widetilde F_T(\mathbf x_i^T)$. The stochastic subfigure shows means over $12$ paired seeds; shaded bands are $95\%$ confidence intervals.}
		\label{fig:hard-instance}
	\end{figure}

	\paragraph{Lower-bound instance}
	We first consider the instance of Harvey et al.~\cite[Section~5]{harvey2019tight} under a deterministic setting. Let $d=T$, $\Omega_T=\{\m z\in\mathbb R^d:\|\m z\|\le1\}$, and $a_j=1/[2(d+1-j)]$. All nodes have the same objective
	\begin{align*}
		F_T(\m z)
		=\frac12\|\m z\|^2+\max_{1\le j\le d+1}\phi_j(\m z),
	\end{align*}
where $\phi_j(\m z)
=\sum_{\ell<j}a_\ell z_\ell-z_j$ for $j\le d$ and
$\phi_{d+1}(\m z)=\sum_{\ell=1}^d a_\ell z_\ell$.
	At $\m z$, the oracle returns
	$
	\m g(\m z)=\m z+\nabla\phi_{j_0(\m z)}(\m z)$, where 
	$j_0(\m z)=\min\arg\max_{1\le j\le d+1}\phi_j(\m z).
	$
	Because the objective, initialization, and oracle are identical across nodes, $\m x_i^t=\bar{\m x}^t$ for every $i$ and $t$; thus this experiment isolates the last-iterate effect. We compare the standard stepsize $\alpha^t=1/t$ with its cooldown version $\alpha^t=p_\ell/t$. Fig.~\ref{fig:hard-instance}a reports the standard and cooldown last iterates and, for the standard schedule, the half-suffix average
	$\widehat{\m x}_{1/2}^T
	=\frac{1}{T-\lfloor T/2\rfloor}
	\sum_{t=\lfloor T/2\rfloor+1}^{T}\bar{\m x}^t$.
	
	For a stochastic variant, augment the variable to $(\m u,\m v)\in\mathbb R^{T+8}$, and use $\widetilde\Omega_T=\{(\m u,\m v):\|\m u\|^2+\|\m v\|^2\le1\}$
	and the objective is 
	$
	\widetilde F_T(\m u,\m v)
	=F_T(\m u)+\|\m v\|^2 / 2+\beta\|\m v\|_1
	$ with $\beta=0.05.$
	The auxiliary term is nonnegative and vanishes at $\m v=\m0$, so $\widetilde F_T^*=F_T^*$. At node $i$, the oracle returns
	$
	\widetilde{\m g}_i^t
	=\left(\m g(\m u_i^t),
	\m v_i^t+\beta\operatorname{sign}(\m v_i^t)+\g\xi_i^t\right),$
    where $
	\g\xi_i^t=0.2\,\varepsilon_i^t\m e_{J_i^t}
	$
	and $\varepsilon_i^t$ is Rademacher, $J_i^t$ is uniform on $\{1,\ldots,8\}$, all samples are independent across nodes and iterations, and $\operatorname{sign}(0)=0$. Hence $\mathbb E[\widetilde{\m g}_i^t\mid\mathcal F_t]\in\partial\widetilde F_T(\m x_i^t)$ and the oracle is bounded on $\widetilde\Omega_T$. The stochastic perturbation creates disagreement across nodes while retaining the lower-bound oracle on the $\m u$-block. 
	
	Fig.~\ref{fig:hard-instance} plots scaled objective gap for $T\in\{2^7,\ldots,2^{11}\}$. In Fig.~1a, the standard last-iterate error increases approximately linearly in $\log T$, whereas the suffix-average and cooldown errors remain nearly constant, consistent with the respective $\mathcal O(\log T/T)$ and $\mathcal O(1/T)$ rates. Fig.~1b shows that the standard network-average and worst-local objective gaps increase with $T$, while their cooldown counterparts remain nearly flat. Under the standard stepsizes, the visible separation between the network-average and worst-local errors reflects terminal network disagreement. Under cooldown schedules, the corresponding separation can be controlled to be smaller; hence the two curves are indistinguishable at the scale of the plot.
	
	\paragraph{Projected hinge-loss problem}
	Each node $i$ holds $m=80$ samples $(\m a_{i,r},b_{i,r})\in\mathbb R^{24}\times\{-1,1\}$. To create heterogeneity, we draw $\mathbb P(b_{i,r}=1)=0.05+0.9{(i-1)}/{(n-1)}$
	and construct features from a fixed sparse direction $\m w_*$ and a fixed drift direction $\m d_*$:
	\[
	\widetilde{\m a}_{i,r}
	=1.2\,b_{i,r}\m w_*+0.85\,\g\xi_{i,r}
	+0.50\frac{i-(n+1)/2}{n-1}\m d_*,
	\qquad
	\g\xi_{i,r}\sim\mathcal N(\m0,\m I_{24}),
	\]
	followed by $\m a_{i,r}=\widetilde{\m a}_{i,r}/\max\{\|\widetilde{\m a}_{i,r}\|,1\}$. On $\Omega=\{\m z:\|\m z\|\le3\}$, define
	$
	f_{i,\mu}(\m z)
	=\frac1m\sum_{r=1}^m[1-b_{i,r}\m a_{i,r}^{\sf T}\m z]_+
	+\frac\mu2\|\m z\|^2,
	$ and
    $
    F_\mu=\frac1n\sum_{i=1}^n f_{i,\mu}.
	$
	At iteration $t$, node $i$ samples $r$ uniformly from $\{1,\ldots,m\}$ and uses
	\[
	\m g_i^t=
	\begin{cases}
		-b_{i,r}\m a_{i,r}+\mu\m x_i^t,
		& b_{i,r}\m a_{i,r}^{\sf T}\m x_i^t<1,\\[1mm]
		\mu\m x_i^t,
		& b_{i,r}\m a_{i,r}^{\sf T}\m x_i^t\ge1.
	\end{cases}
	\]
	Uniform sampling gives $\mathbb E[\m g_i^t\mid\mathcal F_t]\in\partial f_{i,\mu}(\m x_i^t)$, and the normalization and projection imply $\|\m g_i^t\|\le1+3\mu$. We use $(\mu,\gamma^t)=(0,0.65/\sqrt t)$ and $(\mu,\gamma^t)=(0.02,50/(t+80))$
	for the convex and strongly convex cases, respectively, and set $T=3\times10^5$. Centralized projected SGD is run with the same base stepsize and horizon $T$.
	
\begin{figure}[htbp]
		\centering
		\includegraphics[width=.8\linewidth]{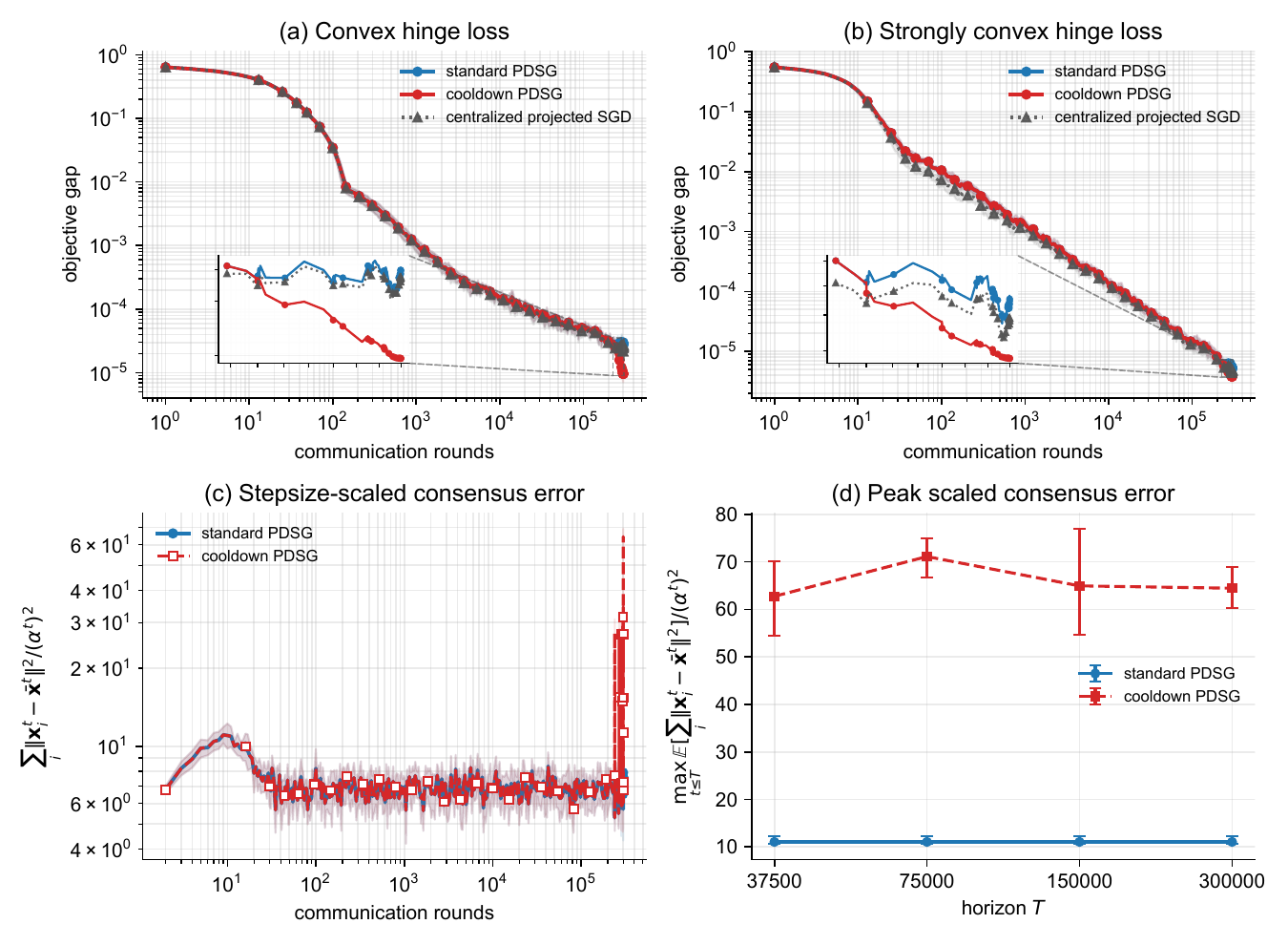}
		\caption{Projected hinge-loss experiments with $8$ paired seeds; shaded bands are $95\%$ confidence intervals. \textnormal{(a)} Convex and \textnormal{(b)} strongly convex objective gaps: the PDSG curves show $n^{-1}\sum_i(F_\mu(\mathbf x_i^t)-F_\mu^*)$, and the centralized curve shows $F_\mu(\mathbf z^t)-F_\mu^*$. \textnormal{(c)} Stepsize-scaled consensus error $\sum_i\|\mathbf x_i^t-\bar{\mathbf x}^t\|^2/(\alpha^t)^2$ in the strongly convex experiment. \textnormal{(d)} Horizon dependence of $\max_{t\le T}\mathbb E[\sum_i\|\mathbf x_i^t-\bar{\mathbf x}^t\|^2]/(\alpha^t)^2$ for $T\in\{37500,75000,150000,300000\}$; the expectation is estimated by the seed average at each $t$ before maximizing over $t\le T$. 
        }
		\label{fig:projected-svm}
	\end{figure}

	Figure~\ref{fig:projected-svm} compares the objective gaps and scaled consensus error, defined precisely in its caption. Each iteration contains one communication round; hence the horizontal axis \emph{communication rounds} is also the iteration number. 
	The three methods have similar objective-gap trajectories until the iteration tail. In the displayed runs, cooldown PDSG reduces terminal oscillation and yields final gaps comparable to or smaller than those of standard PDSG. 
	At the cooldown boundary, the coefficient in $\alpha^t$ halves. Because the consensus $\sum_i\|\mathbf x_i^t-\bar{\mathbf x}^t\|^2$ is evaluated before the $t$-th iterative update, it still reflects the preceding iterate and changes little at the boundary; the scaled error therefore rises. 
	Across the four horizons in Figure~\ref{fig:projected-svm}d, the peak scaled consensus error stays between $62.7$ and $71.1$ for cooldown PDSG, supporting a horizon-uniform constant.
	
	\section{Conclusion}\label{sec:conclusion}
This paper clarifies the last-iterate behavior of PDSG in decentralized
nonsmooth stochastic optimization. Under standard stepsizes, harmonic
accumulations in the descent and last-iterate estimates lead to the
$\mathcal{O}(\log T/\sqrt{T})$ and
$\mathcal{O}(\log T/T)$ rates. In decentralized settings, nonsmoothness
and data heterogeneity may additionally sustain nonvanishing local
directions near a global solution, making the terminal iterates
particularly sensitive to the stepsize schedule.
Our main result shows that removing this logarithmic loss does not
require a particular dyadic construction. We develop a general
multistage cooldown framework in which the stepsize variation is
compatible with network contraction and the objective gap can be
propagated across successive segments. This yields the order-optimal
$\mathcal{O}(1/\sqrt{T})$ and $\mathcal{O}(1/T)$ rates for every last
local iterate and provides a general stepsize-design principle for
decentralized stochastic subgradient methods.
 Numerical experiments show that the cooldown schedule reduces terminal oscillations and improves the terminal accuracy while maintaining controlled consensus error.

\section*{Acknowledgments}
We thank Prof. Hongchao Zhang from Louisiana State University for  constructive suggestions during the preparation of this work.
	
	\appendix

	\section{Key lemmas}
	\begin{lemma}\label{lem:mod-initial-dist}
		Suppose $F$ is $\mu$-strongly convex \newredfix{on $\Omega$} and the base stepsize is $\gamma^t=2/(\mu t)$. During the initial phase $1\le t\le T_1$, for which $\alpha^t=\gamma^t$, define $D^t:=n^{-1}\mathbb{E}[\|\m{x}^t-\m{1}\otimes\m{z}^*\|^2]$.
		Then, we have $
			D^t\le {4Q_{sc}^{mod}}/{\mu^2t},\quad 1\le t\le T_1.
		$
	\end{lemma}
	
	
	\begin{lemma}\label{lem_base_sc_optimal}
		Suppose the global objective $F$ is $\mu$-strongly convex \newredfix{on $\Omega$} and the base stepsize is $\gamma^t = {2}/{(\mu t)}$. Then the iterate indexed by $\tau_0$, which minimizes the expected objective over the suffix window $[\lfloor T_1/2 \rfloor + 1, T_1]$, satisfies the following order-optimal $\mathcal{O}(1/T_1)$ bound:
		\begin{align*}
			\mathbb{E}[F(\bar{\m{x}}^{\tau_0})] - F^* \newredfix{\leq \frac{2 Q_{sc}^{mod}(1+\ln 2) + 4K_{mod}}{\mu T_1}} = \mathcal{O}\left(\frac{1}{T_1}\right).
		\end{align*}
	\end{lemma}
	
	\begin{proof}
		The proof uses the half-suffix averaging argument of Rakhlin et al. \cite{rakhlin2012making}. Because $\tau_0$ minimizes the expected objective over this suffix window, its objective gap is at most the averaged gap:
		\begin{align}\label{proof_sc_inf_avg}
			&\mathbb{E}[F(\bar{\m{x}}^{\tau_0})] - F^* = \min_{\lfloor T_1/2 \rfloor < t \leq T_1} \mathbb{E}[F(\bar{\m{x}}^t)- F^*] \\
			\leq &\frac{1}{T_1-\lfloor T_1/2 \rfloor} \sum_{t=\lfloor T_1/2 \rfloor + 1}^{T_1} \mathbb{E}[F(\bar{\m{x}}^t)- F^*].\notag
		\end{align}
		Summing the inequality \eqref{Fxt-Fu-mod} over  $t \in [\lfloor T_1/2 \rfloor + 1, T_1]$ yields
		\begin{align}\label{lemB_eq2}
			&\sum_{t=\lfloor T_1/2 \rfloor + 1}^{T_1} \mathbb{E}[F(\bar{\m{x}}^t) - F^*] \leq \frac{D^{\lfloor T_1/2 \rfloor + 1}}{2\gamma^{\lfloor T_1/2 \rfloor + 1}} \\
			&+ \sum_{t=\lfloor T_1/2 \rfloor + 2}^{T_1} \frac{D^t}{2} \left( \frac{1}{\gamma^t} - \frac{1}{\gamma^{t-1}} \right) + K_{mod} \sum_{t=\lfloor T_1/2 \rfloor + 1}^{T_1} \gamma^t.\notag
		\end{align}
		Substituting $\gamma^t = {2}/{(\mu t)}$, for which ${1}/{\gamma^t} - {1}/{\gamma^{t-1}} = {\mu}/{2}$, and applying Lemma~\ref{lem:mod-initial-dist} bounds the RHS of \eqref{lemB_eq2} as follows:
		\begin{align*}
			\text{RHS} &\leq \frac{Q_{sc}^{mod}}{\mu} + \sum_{t=\lfloor T_1/2 \rfloor + 2}^{T_1} \frac{Q_{sc}^{mod}}{\mu t} + \sum_{t=\lfloor T_1/2 \rfloor + 1}^{T_1} \frac{2K_{mod}}{\mu t}\\
			& \newredfix{= \frac{Q_{sc}^{mod}}{\mu}\left(1+\sum_{t=\lfloor T_1/2 \rfloor + 2}^{T_1}\frac{1}{t}\right) + \frac{2K_{mod}}{\mu}\sum_{t=\lfloor T_1/2 \rfloor + 1}^{T_1}\frac{1}{t}.}
		\end{align*}
		\newredfix{Applying the following inequalities that $\sum_{t=\lfloor T_1/2 \rfloor + 1}^{T_1} \frac{1}{t} \leq \frac{(T_1-\lfloor T_1/2 \rfloor)}{(\lfloor T_1/2 \rfloor+1)} \leq 1$ and}
		\begin{align*}
			\sum_{t=\lfloor T_1/2 \rfloor + 2}^{T_1} \frac{1}{t} \leq \int_{\lfloor T_1/2 \rfloor+1}^{T_1} \frac{1}{x} \text{d}x = \ln\left( \frac{T_1}{\lfloor T_1/2 \rfloor+1} \right) \leq \ln 2,
		\end{align*}
		\newredfix{we obtain}
		\newredfix{$
			\sum_{t=\lfloor T_1/2 \rfloor + 1}^{T_1} \mathbb{E}[F(\bar{\m{x}}^t) - F^*] \leq (Q_{sc}^{mod}(1+\ln 2)+2K_{mod})/{\mu}.
		$}
		Dividing by $T_1 - \lfloor T_1/2 \rfloor \geq {T_1}/{2}$ and using the result in \eqref{proof_sc_inf_avg} completes the proof.
	\end{proof}
	
	\section{Proof of Lemma \ref{lem_consensus}}\label{proof_lem_consensus}
	
	Set \(a:=s+q\) and
	$
	\Gamma_a(\lambda):=\sum_{k=0}^{\infty}(1+k)^{2a}\lambda^k<\infty.
	$
	Assumption~\ref{as_mix}
	and \(\gamma^j/\gamma^t=(t/j)^q\le(1+t-j)^q\) give
	\begin{equation}\label{eq_stepsize_ratio_new}
		\left(\frac{\alpha^j}{\alpha^t}\right)^2
		\le C_p^2(1+t-j)^{2a},
		\qquad 1\le j\le t\le T.
	\end{equation}
	Taking \(j=1\) in Assumption~\ref{as_mix} yields
	\(p_{i(t)}\ge C_p^{-1}t^{-s}\), and hence
	\(\alpha^t\ge cC_p^{-1}t^{-a}\).
	Let
	\(Q^{t+1}:=\mathbb{E}[\|\m{x}^{t+1}-\m{M}\m{x}^{t+1}\|^2]\).
	\newredfix{For \(0<\sigma<1\), set $\lambda:=\frac{1+\sigma^2}{2}$, $\bar C:=\frac{1+\sigma^2}{1-\sigma^2}nG^2$,
	which follow by applying Young's inequality with parameter
	\((1-\sigma^2)/(2\sigma^2)\). If \(\sigma=0\), then
	\(\m W=\m M\), giving $Q^{t+1}\le(\alpha^t)^2\mathbb{E}[\|(\m I-\m M)\m g^t\|^2]
	\le nG^2(\alpha^t)^2$.
	In this case, set \(\lambda:=1/2\) and \(\bar C:=nG^2\). Thus, in either case,
	\(0<\lambda<1\) and
	\(Q^{t+1}\le\lambda Q^t+\bar C(\alpha^t)^2\).}
	Therefore
	\begin{equation}\label{lem_decay_aux_eq1}
		\lambda^t\le C_1(\alpha^t)^2,
		\qquad
		C_1:={C_p^2}/{c^2}\sup_{r\ge1}r^{2a}\lambda^r<\infty.
	\end{equation}
	Unrolling this recursion and applying \eqref{eq_stepsize_ratio_new} yield
	$
		Q^{t+1}
		\le\lambda^tQ^1+
		\bar C\sum_{j=1}^t(\alpha^j)^2\lambda^{t-j}
		\le\lambda^tQ^1+
		\bar C C_p^2\Gamma_a(\lambda)(\alpha^t)^2.
	$
	Together with \eqref{lem_decay_aux_eq1}, this gives
	\begin{equation}\label{Qt+1_bound_v2}
		Q^{t+1}\le\Cmodraw(\alpha^t)^2,
		\qquad
		\Cmodraw:=C_1Q^1+\bar C C_p^2\Gamma_a(\lambda).
	\end{equation}
	Finally, Assumption~\ref{as3} and \(\gamma^t=ct^{-q}\) imply
	\(\alpha^{t-1}\le c_\alpha\alpha^t\) for \(t\ge2\), where one may
	take \(1\le c_\alpha\le2^q/\rho\). Applying
	\eqref{Qt+1_bound_v2} at iteration \(t-1\) gives
	$
	\mathbb{E}[\|\m{x}_i^t-\bar{\m{x}}^t\|^2]
	\le Q^t\le\Cmod(\alpha^t)^2,
	$
	with $\Cmod:=\max\left\{\Cmodraw c_\alpha^2,
	{Q^1}/{(\alpha^1)^2}\right\}$. This proves the claim.
	
	\bibliographystyle{siamplain}
	\bibliography{refs}
\end{document}